\documentclass[reqno,11pt]{amsart}

\usepackage{amssymb, amsmath,amsthm,amscd}
\usepackage{mathtools}

\usepackage{natbib}
\usepackage{cite}
\usepackage{enumitem}
\usepackage[
  margin=2.8cm,
  includefoot,
  footskip=30pt,
]{geometry}
\usepackage{ytableau}
\usepackage{graphicx}

\numberwithin{equation}{section}

\newtheorem{Thm}{Theorem}[section]
\newtheorem{Lemma}[Thm]{Lemma}
\newtheorem{Prop}[Thm]{Proposition}
\newtheorem{Cor}[Thm]{Corollary}

\theoremstyle{definition}
\newtheorem{Def}[Thm]{Definition}
\newtheorem{Ex}[Thm]{Example}

\theoremstyle{remark}

\newtheorem{Remark}[Thm]{Remark}

\usepackage{hyperref}

\makeatletter
\newcommand\incircbin
{%
  \mathpalette\@incircbin
}
\newcommand\@incircbin[2]
{%
  \mathbin%
  {%
    \ooalign{\hidewidth$#1#2$\hidewidth\crcr$#1\bigcirc$}%
  }%
}

\makeatother

\newcommand{\pref}[1]{(\ref{#1})}

\DeclarePairedDelimiter\floor{\lfloor}{\rfloor}

\newcommand{\N}{\mathbb{N}}
\newcommand{\R}{\mathbb{R}}
\newcommand{\C}{\mathbb{C}}

\renewcommand{\SS}{\mathbb{S}}

\newcommand{\mg}{\mathfrak{g}}

\renewcommand{\sl}{\mathfrak{sl}}
\newcommand{\so}{\mathfrak{so}}

\newcommand{\OLie}{\mathrm{O}}
\newcommand{\SO}{\mathrm{SO}}
\newcommand{\SL}{\mathrm{SL}}

\newcommand{\local}[2]{\left. {#1} \right|_{#2}}

\newcommand{\p}{\partial}

\newcommand{\pd}[2]{\frac{\p{#1}}{\p{#2}}}

\newcommand{\pdxi}[2]{\pd{#1}{x^{#2}}}

\newcommand{\df}[1]{d{#1}}
\newcommand{\dxi}[1]{\df{x^{#1}}}
\newcommand{\dyi}[1]{\df{y^{#1}}}

\newcommand{\Char}[1]{\chi_{\raisebox{-1mm}{\scriptsize ${#1}$}}}

\DeclareMathOperator{\id}{id}

\DeclareMathOperator{\sgn}{sgn}

\DeclareMathOperator{\Hom}{Hom}

\newcommand{\skprod}[1]{\left\langle{#1}\right\rangle}

\newcommand{\KK}{\mathcal{K}}

\newcommand{\BB}{\mathcal{B}}

\newcommand{\Lef}{L}
\newcommand{\DLef}{\Lambda}

\newcommand{\AMF}[1]{\textstyle{\bigwedge^{#1}}}

\newcommand{\AMFP}[1]{\textstyle{\bigwedge_p^{#1}}}

\newcommand{\SDF}[3]{\Omega^{#1}({#2})^{#3}}
\newcommand{\SDFH}[3]{\Omega_h^{#1}({#2})^{#3}}
\newcommand{\SDFV}[3]{\Omega_v^{#1}({#2})^{#3}}

\DeclareMathOperator{\nc}{nc}

\DeclareMathOperator{\Val}{Val}
\DeclareMathOperator{\TVal}{Val}
\newcommand{\Valsm}{\Val^{sm}}

\DeclareMathOperator{\Curv}{Curv}
\DeclareMathOperator{\TCurv}{Curv}

\newcommand{\TInv}[1]{\overline{#1}}

\DeclareMathOperator{\integ}{integ}
\DeclareMathOperator{\glob}{glob}
\DeclareMathOperator{\Integ}{Integ}

\DeclareMathOperator{\vol}{vol}

\DeclareMathOperator{\Sym}{Sym}

\newcommand{\SymL}[1]{\Sym^2\!\Lambda^{{#1}}}
\DeclareMathOperator{\Ind}{Ind}
\DeclareMathOperator{\Res}{Res}

\DeclareMathOperator{\tr}{tr}

\newcommand{\FPhi}{\tilde{\Phi}}
\newcommand{\FPsi}{\tilde{\Psi}}
\newcommand{\FXi}{\tilde{\Xi}}
\newcommand{\FTheta}{\tilde{\Theta}}

\newcommand{\FT}{\tilde{T}}

\newcommand{\CPhi}{\Phi}
\newcommand{\CPsi}{\Psi}
\newcommand{\CXi}{\Xi}
\newcommand{\CTheta}{\Theta}
\newcommand{\CT}{T}
\newcommand{\VPhi}{\phi}
\newcommand{\VPsi}{\psi}

\newcommand{\VT}{\tau}

\usepackage{xcolor}

\begin{document}
\bibliographystyle{plain}
\setcitestyle{numbers,open={[},close={]},citesep={,}} 

\begin{abstract}
Valuations constitute a class of functionals on convex bodies which include the Euler-characteristic, the surface area, the Lebesgue-measure, and many more classical functionals. Curvature measures may be regarded as "localised`` versions of valuations which yield local information about the geometry of a body's boundary. 

A complete classification of continuous translation-invariant $\SO(n)$-invariant valuations and curvature measures with values in $\R$ was obtained by Hadwiger and Schneider, respectively. More recently, characterisation results have been achieved for curvature measures with values in $\Sym^p \R^n$ and $\SymL{q} \R^n$ for $p,q \geq 1$ with varying assumptions as for their invariance properties. 

In the present work, we classify all smooth translation-invariant $\SO(n)$-covariant curvature measures with values in any $\SO(n)$-representation in terms of certain differential forms on the sphere bundle $S\R^n$ and describe their behaviour under the globalisation map. The latter result also yields a similar classification of all continuous $\SO(n)$-covariant valuations with values in any $\SO(n)$-representation. Furthermore, a decomposition of the space of smooth translation-invariant $\R$-valued curvature measures as an $\SO(n)$-representation is obtained. As a corollary, we construct an explicit basis of continuous translation-invariant $\R$-valued valuations.
\end{abstract}

\title{Characterisation of Valuations and Curvature Measures in Euclidean Spaces}
\author{Mykhailo Saienko}
\email{saienko@math.uni-frankfurt.de}

\address{Institut f\"ur Mathematik, Goethe-Universit\"at Frankfurt,
Robert-Mayer-Str. 10, 60054 Frankfurt, Germany}
\thanks{Supported by DFG grants BE 2484/5-1 and BE 2484/5-2.\\ 
AMS 2010 
}
\maketitle

\section{Introduction}
\setcounter{page}{1}

\subsection{Background}
\label{sec:intro-background}
Let $\KK(\R^n)$ be the set of convex bodies, i.e., compact convex subsets, in $\R^n$ and $A$ be an Abelian semigroup. The map $\phi : \KK(\R^n) \to A$ is called a \emph{valuation} if it satisfies the equation:
\begin{equation}
	\phi(K\cup L) + \phi(K\cap L)= \phi(K)+\phi(L).
	\label{eq:IncExclPrinciple}
\end{equation}
whenever $K\cup L \in \mathcal{K}(\R^n)$. We study the case where $A$ is a finite-dimensional Euclidean vector space $\R^m$. A valuation $\phi$ is then said to be continuous if it is continuous with respect to the topology induced by the Hausdorff-metric on $\KK(\R^n)$. Valuations can be studied on broader classes of subsets in $\R^n$ or on certain subsets of manifolds \citep{Alesker2006I, Alesker2006II, Bernig2012, BernigLefschetz, FuKinematicWDC}. Other important target spaces include the case $A=\KK(\R^n)$ (Minkowski valuations) \citep{LudwigMinkowski} and the space of signed measures on the sphere (area measures) \citep{WannererAreaIntGeo, WannererAreaModule}.

The first valuations to become objects of systematic study were continuous $\R$-valued valuations invariant under the action of the Euclidean group $\TInv{\SO(n)}:= \SO(n)\ltimes \R^n$. Hadwiger \citep{HadwigerIntGeo} showed them to form an $(n+1)$-dimensional vector space $\Val^{\SO(n)}$ spanned by the \emph{intrinsic volumes} $\mu_0,\ldots,\mu_n$, where $\mu_0$ is the Euler-characteristic and $\mu_n$ is the Lebesgue-measure. 

Almost 50 years later, Alesker initiated the program of describing continuous valuations invariant -- but also equi- and contravariant -- under different Lie-groups $G$. It resulted in a range of Hadwiger-type results \citep{AbardiaDiffBodies2012, Alesker2003, Alesker2008, BernigGAFA, BernigSOpq, HIG, BernigSolanesQuaternions, LudwigReitzner99, LudwigReitzner10, WannererContravariant, Solanes15,  WannererEquivariant, WannererAreaIntGeo, WannererAreaModule}.

Dropping $G$-invariance, the space $\Val$ of continuous translation-invariant valuations was shown by McMullen in \citep{McMullen} to admit a decomposition by homogeneity degree and parity:
$$ \Val = \bigoplus_{0\leq k\leq n} \Val^{+}_k \oplus \Val^-_k,$$
where $\Val^{\pm}_k$ are \emph{infinite-dimensional} (Fr\'echet-)spaces unless $k\in \{0,n\}$, in which case $\Val_k$ is one-dimensional and spanned by the Euler-characteristic and the Lebesgue-measure, respectively. 

A more advanced structure result is the decomposition of $\Val_k$ in $\SO(n)$-irreducible representations by Alesker, Bernig, and Schuster \citep{AleskerBernigSchuster}. They showed that $\Val_k$ is \emph{multiplicity-free} and contains the irreducible $\SO(n)$-representations $\Gamma_{[\lambda]}$ with highest weights $\lambda$ such that:
\begin{itemize}
	\item $\lambda_j = 0$ for $j > \min(k,n-k)$;
	\item $|\lambda_j| \neq 1$ for $1\leq j \leq \floor{n/2}$;
	\item $|\lambda_2| \leq 2$.
\end{itemize}
The parity of the valuation corresponds to the parity of $\lambda_1$, while the case $\lambda_2=0$ corresponds to the so-called spherical valuations. The $\Gamma^*$-typical component in this decomposition can be identified with the space 
\begin{equation}
	\label{eq:ValkGamma}
\TVal_{k,\Gamma}^{\SO(n)}:=\{\phi \in \Val_k \otimes \Gamma: \phi(K)=g \phi(g^{-1}K), K \in \KK(\R^n), g \in \SO(n)\}. 
\end{equation}

The explicit bases for $\TVal_{k,\Gamma}^{\SO(n)}$ have remained elusive for several years. In fact, the structure of $\TVal_{k,\Gamma}^{\SO(n)}$ is only known for $\Gamma=\Sym^p \R^n$, as several bases and global kinematic formulae were gradually elaborated by different authors, including Alesker, Bernig, Hug, McMullen, and Schuster \citep{AleskerTensor99, BernigHugTensor2013, HugSchneiderSchuster2007, HadwigerSchneider71, Hug2018KineticSymVal, HugSchneiderSchuster2008, McMullenTensor}.

The present paper closes this gap by establishing in rather explicit terms a basis of $\TVal_{k,\Gamma}^{\SO(n)}$ for \textit{any} $\SO(n)$-representation $\Gamma$. To achieve this, we extend our study to \textit{curvature measures}, an extremely useful concept through which the study of continuous translation-invariant valuations can be linked to the more familiar concepts of differential forms on the sphere bundle $S\R^n$. Let us briefly outline this connection.

Curvature measures were introduced by Federer in an attempt to connect several integral-geometric results that had been previously disparate \citep{FedererCurv}. He observed that intrinsic volumes $\mu_k(K), k=0,\ldots,n-1$ can be computed by integrating the symmetric functions of the principal curvatures over its boundary $\partial K$ if it is sufficiently smooth. Replacing $\partial K$ under the integral by $\partial K \cap U$ for any Borel-set $U$, one naturally obtains a ``localised" version of $\mu_k$ called the $k$-th \emph{Lipschitz-Killing curvature measure} $\Phi_k : \KK(\R^n)\times \BB(\R^n) \to \R$, where $\BB(\R^n)$ is the Borel-$\sigma$-Algebra on $\R^n$.  Obviously, $\mu_k$ can be recovered from $\Phi_k$ by the relation $\mu_k(K) = \Phi_k(K, \R^n)$ for any $K\in\KK(\R^n)$. It is by no means trivial to extend this description of Lipschitz-Killing curvature measures to non-smooth convex bodies. In fact, this was one of the main results of Federer's publication and a major driving force to developing the geometric measure theory.

The name ``curvature measures" is more than justified for $\Phi_k$. On the one hand, for $K$ sufficiently smooth and $U$ any Borel-set, $\Phi_k(K,U)$ yields local information about the curvature of $\partial K$. On the other hand, $\Phi_k(K, \cdot)$ is a non-negative Borel-measure for a fixed convex body $K$ that is \textit{weakly-continuous}, i.e.:
$$ \int_{\R^n} f(x) d\Phi_k(K_i, x) \to \int_{\R^n} f(x) d\Phi_k(K, x)$$
for any continuous function $f:\R^n \to \R$ and any sequence of convex bodies $K_i$ converging to a convex body $K$ \citep[pp. 288ff.]{Schneider}. The ``localisation" procedure also preserves the $\TInv{\SO(n)}$-invariance of $\Phi_k$, i.e., $\Phi_k(\TInv{g}K, \TInv{g}U) = \Phi_k(K,U)$ for all $\TInv{g}\in\TInv{\SO(n)}$, $K\in\KK(R^n)$, $U\in\BB(\R^n)$. In fact, $\Phi_k$ comprise the basis of $\TInv{\SO(n)}$-invariant weakly continuous curvature measures $\Curv^{\SO(n)}$ on convex bodies in $\R^n$ \citep{SchneiderCurv}. 

Later, Z\"ahle \citep{Zaehle86} discovered that $\Phi_k$ and $\mu_k$ for all $k<n$ can be represented as
\begin{equation}\label{eq:normal-cycle}
\Phi_k(K,U)=\int_{\nc(K) \cap \pi^{-1}(U)} \omega_k\quad\hbox{ and } \quad \mu_k(K) = \int_{\nc(K)} \omega_k, 
\end{equation}
where $\nc(K)$ is a Lipschitz-submanifold of the sphere bundle $S\R^n$ called the \textit{normal cycle} of $K$, $\pi: S\R^n \to \R^n$ is the natural projection and $\omega_k$ is a certain $\TInv{\SO(n)}$-invariant differential form on $S\R^n$ of bi-degree $(k, n-k-1)$. Replacing $\omega_k$ with any translation-invariant form $\omega \in \Omega^{n-1}(S\R^n)$, the functional
$$\Phi_\omega(K,U) := \int_{\nc(K) \cap \pi^{-1}(U)} \omega$$
induces a continuous translation-invariant valuation $K\mapsto \Phi_\omega(K,\R^n)$ and a weakly continuous translation-invariant Borel-measure $(K,U)\mapsto \Phi_\omega(K,U)$. The former are called smooth valuations and the vector space spanned by them is denoted by $\Valsm$. The latter are referred to as \textit{smooth translation-invariant curvature measures}. We will denote the vector space formed by them by $\Curv^{sm}$. The valuation $\phi_\omega(\cdot) := \glob(\Phi_\omega)(\cdot):=\Phi_\omega(\cdot,\R^n)$ is called the \textit{globalisation} of $\Phi_\omega$ and the globalisation map $\glob : \Curv^{sm} \to \Valsm$ is trivially onto. However, contrary to the case of $\mu_k$ and $\Phi_k$, the kernel of $\glob$ is not trivial, i.e., the ``localisation" procedure is not canonical.

The space $\Valsm$ possesses rich algebraic structures, such as product, convolution and a Fourier-type transform \citep{Alesker2011Fourier, Alesker2006III, Bernig2006}, that are connected to the kinematic formulae \citep{FuFTAIG} and allow to write out such formulas explicitly \citep{BernigGAFA, bernig_g2, HIG, bernig_solanes17}. Furthermore, a corollary of Alesker's famous Irreducibility Theorem \citep{Alesker2001} states that $\Valsm$ lies densely in $\Val$ and, in particular, that all valuations from the finite-dimensional space $\TVal_{k,\Gamma}^{\SO(n)}$ are smooth (Proposition \ref{thm:smoothness}). 

It is this fact and the careful examination of the kernel of the globalisation map (Theorem \ref{thm:BernigLemma}) that allow us to describe the basis of $\TVal_{k,\Gamma}^{\SO(n)}$ in terms of the basis of the space $\TCurv_{k,\Gamma}^{sm, \SO(n)}$ of smooth $\SO(n)$-covariant translation-invariant curvature measures with values in $\Gamma$  (Proposition \ref{thm:TValBasis}). Establishing the latter is the main result of this work  (Theorem \ref{thm:OnHarmonicBasis}) and requires, among other mathematical tools, the harmonic decomposition of $\Curv^{sm}$  (Theorem \ref{thm:curvHarmonicDecomposition}).

Our work (Remark \ref{rem:ThetaSObutNotO}) has revealed that -- surprisingly and in contrast to $\SO(n)$-equivariant translation-invariant valuations -- there are $\SO(n)$-equivariant curvature measures that are not $\OLie(n)$-equivariant. This has entailed new efforts to classify them for $\Gamma=\Sym^p\R^n$ on convex polytopes and to study their extensions to convex bodies \citep{HugSchneiderTensorCurvs2017, HugSchneiderTensorCurvs2016}. Furthermore, we show in Proposition \ref{thm:formsGeometric} that the differential forms constructed in our work are intimately related to those used to classify the so-called local Minkowski-tensors with certain properties \citep{HugSchneiderTensorCurvs2013} and later to establish several integral-geometric formulae for them \citep{Hug2019KineticSymCurv}. 

Finally, we complete the search for smooth $\SO(n)$-covariant translation-invariant curvature measures with values $\Gamma = \SymL{q} \R^n$ started by Bernig in \citep{BernigTensor2005} discover more symmetries for them (Proposition \ref{thm:BernigCurvs} and Proposition \ref{thm:exchangeColumns}).

The bases of $\TVal_{k,\Gamma}^{\SO(n)}$ also induce a \textit{Schauder}-basis of $\Val$ (Proposition \ref{thm:TValBasis}). This might turn useful for a range of applications. For example, a famous result by Klain \citep{KlainRota} states that $\Val_k^+$ can be seen as a subspace of the space of functions on the Grassmannian of $k$-planes in $V$. This allows to relate the basic operators on $\Valsm$ -- such as the Lefschetz operator, i.e., multiplication by the first intrinsic volume, and the derivation operator, i.e., convolution with the $(n-1)$-st intrinsic volume -- to some known integral transforms on Grassmannians, such as the Radon transform and the cosine transform. This approach has lead to some deep results \citep{alesker_hodge_riemann, Alesker2003, Alesker2004, BernigHugTensor2013, BernigSolanesQuaternions, bernig_solanes17, dorrek_schuster, KotrbatyHodgeRiemann, kotrbaty_wannerer} in the \textit{even} case. 

These results cannot be easily extended to the \textit{odd} case, as there is no embedding for odd valuations which would be comparable to Klain's map. Bernig and Hug studied in \citep{BernigHugTensor2013} \emph{spherical valuations} and proved kinematic formulas for tensor valuations using tools from harmonic analysis. Although spherical valuations may be of odd parity, they do not form a dense subspace in $\Val$. Our hope is that the basis of $\Val$ we discovered -- being compatible with the harmonic decomposition of $\Val$ and thus allowing for very precise control of the parity of its elements -- might serve the same function for the odd case as Klain's map did for even valuations.

The plan of the paper is as follows. In Subsection \ref{sec:intro-results}, we formulate the main results of this work. In Section \ref{sec:RepTheory}, we recall all necessary basics of the finite-dimensional representation theory of $\SL(n)$ and $\SO(n)$, including Young-symmetrisers, trace-free spaces as well as restricted and induced representations. We refer to \citep{Fulton, FultonHarris, GoodmanWallach} for more detailed expositions of this topic. In Section \ref{sec:valuations} we discuss some facts from the valuation theory which we need to prove the main results. The prominent references here are \citep{BernigAIG, AIG, KlainRota, Schneider} along with the papers mentioned above. The new results are proven in Section \ref{sec:coreResults}.

\subsection{Main Results} \label{sec:intro-results}
The space $\Curv^{sm}_{k}$ naturally admits the structure of an $\SO(n)$-module by $ (g \Phi)(K, U) := \Phi(g^{-1} K, g^{-1} U)$ for all $K\in \KK(\R^n),\ U\in \BB(\R^n)$. By the Theorem of Peter-Weyl, $\Curv^{sm}_k$ may be written as a direct sum of irreducible finite-dimensional $\SO(n)$-modules. All such $\SO(n)$-representations may be uniquely characterised up to isomorphism by tuples $\lambda=(\lambda_1\geq \ldots\geq \lambda_{\floor{n/2}})$ such that $\lambda_{\floor{n/2}}\geq 0$ if $n$ is odd and $\lambda_{n/2-1} \geq |\lambda_{n/2}| \geq 0$ if $n$ is even. 

\begin{Thm}
\label{thm:curvHarmonicDecomposition}
Let $n\geq 2$, $0\leq k \leq n-1$. Then $\Curv^{sm}_k$ consists precisely of $\SO(n)$-representations $\Gamma_{[\lambda]}$ with tuples $\lambda$ such that:
\begin{itemize}
	\item $\lambda_j = 0$ for $j > \min(k+1,n-k)$;
	\item $|\lambda_j| = 1$ for at most one $1\leq j\leq \floor{n/2}$;
	\item $|\lambda_2| \leq 2$.
\end{itemize}
Let $m$ be the highest $j$ such that $\lambda_j\neq 0$. The multiplicity $\Gamma_{[\lambda]}$ in $\Curv^{sm}$ is $2$ except if $m=\min(k+1,n-k)$ or $|\lambda_m| < 2$ (in which case it is $1$) and if $n= 2k + 1, m = k, |\lambda_m| \geq 2$ in which case it is $3$.
\end{Thm}

We now turn to constructing the basis of  $\TCurv_{k,\Gamma}^{sm, \SO(n)}$. Let $e_{i}$, $i=1,\ldots, n$ be the standard orthonormal basis of $\R^n$, $\dxi{i},\dyi{i}$ be the canonical frame on the cotangent bundle $T^*\R^n$ and write $e_{\otimes i_1,\ldots,i_q y} := e_{i_1}\otimes\ldots\otimes e_{i_q}\otimes y$. Define for $0\leq k\leq n-1$, $p \geq 0$ and $0 \leq q \leq \min(k,n-k-1)$ the following families of differential forms pointwise for $(x,y)\in S\R^n$:
\begin{equation}
\begin{aligned}
\FPhi^n_{\otimes k,p,q} & =  C_n \sum_{\pi,i} \sgn{\pi}\ y_{\pi_{n}}\, \dxi{i_1\ldots i_q \pi_{q+1}\ldots \pi_{k}}\wedge\dyi{\pi_{k+1}\ldots\pi_{n-1}} \otimes e_{\otimes i_1\ldots i_{q}}\otimes\,e_{\otimes \pi_{1}\ldots\pi_{q}}\otimes y^{p}, \\
\FXi^n_{\otimes k,p,q} & = C_n  \sum_{\pi,i} \sgn{\pi}\ y_{\pi_{n}}\, \dxi{i_1\ldots i_q \pi_{q+1}\ldots \pi_{k}}\wedge\dyi{\pi_{k+1}\ldots\pi_{n-1}} \otimes e_{\otimes i_1\ldots i_{q}y}\otimes\,e_{\otimes \pi_{1}\ldots\pi_{q}}\otimes y^{p}, \\
\FPsi^n_{\otimes k,p,q+1} & = C_n \sum_{\pi,i} \sgn{\pi}\ y_{\pi_{n}}\, \dxi{i_1\ldots i_{q} \pi_{q+1}\ldots \pi_{k}}\wedge\dyi{\pi_{k+1}\ldots\pi_{n-1}}\otimes e_{\otimes i_1\ldots i_{q}y} \otimes\,e_{\otimes \pi_{1}\ldots\pi_{q}y}\otimes y^{p},
\end{aligned}
\label{eq:FormDefs}
\end{equation}
where $C_n= (-1)^{n-1}$ and we sum over all $n$-permutations $\pi\in S_n$ and indexes $i_1,\ldots i_q=1,\ldots, n$. The above forms assume values in $(\R^n)^{\otimes 2q+p}$, $(\R^n)^{\otimes 2q+p+1}$, and $(\R^n)^{\otimes 2q+p+2}$, respectively. Additionally, define for $k\geq 1$, $n=2k+1$, and $p\geq 0$ a family of $(\R^n)^{\otimes 2k+p}$-valued forms:
\begin{eqnarray*}
\FTheta^n_{\otimes k,p} & = & \sum_{i,j} \dxi{i_1\ldots i_k}\wedge\dyi{j_1\ldots j_k} \otimes e_{\otimes i_1\ldots i_k}\, \otimes\, e_{\otimes j_1\ldots j_k}\, \otimes y^{p},
\end{eqnarray*}
where the sum is over the indexes $i_1,\ldots i_k,j_1,\ldots j_k=1,\ldots,n$. We will often omit the superscript $n$ and use $T$ as a generic letter that may stand for $\CPhi, \CXi, \CPsi,$ or $\CTheta$. 

Special cases of such forms have been used before in different contexts. Write $\CT_{\otimes k,p,q}$ for the curvature measure induced by $\FT_{\otimes k,p,q}$.

\begin{Prop}
\label{thm:BernigCurvs}
Let $\Psi_{k,d}$ be the $\SymL{d}\R^n$-valued curvature measures defined in \citep{BernigTensor2005}. Then:
\begin{equation*}
	\Psi_{k,d} = \frac{1}{s_{n-k-1} (k-d)!\,d!\, (n-k-1)!}\, \CPhi_{\otimes k,0,d}.
\end{equation*}	
\end{Prop}

\begin{Prop} \label{thm:formsGeometric}
Let $P\subset \R^n$ be an arbitrary convex polytope and denote by $\mathcal{F}_k(P)$ the set of all its $k$-dimensional faces.

Let $W \subset \R^n$ be a $k$-dimensional vector subspace and write $Q_W$ for the restriction to $W$ of the metric tensor $Q$ preserved by $\OLie(n)$. Taking $v_1,\ldots, v_k$ to be an orthonormal basis of $W$ so that $Q_W = \sum_{i=1}^k v_i\otimes v_i$ and writing $v_{i_1\ldots i_q}:= v_{i_1}\wedge \ldots \wedge v_{i_q}$, define $Q_W^{\wedge q} := \sum_{i_1,\ldots i_q=1}^k v_{i_1\ldots i_q}\otimes v_{i_1\ldots i_q} \subset \AMF{q}\R^n \otimes \AMF{q}\R^n \subset (\R^n)^{\otimes 2q}$ the $q$-fold wedge product of $Q_W$ with itself. Then,
$$\Phi_{\otimes k,p,q}(P,U) = (-1)^{n-1}\frac{(k-q)!(n-k-1)!}{q!} \sum_{F\in \mathcal{F}_k} \vol (F \cap U)\,Q_{L(F)}^{\wedge q} \otimes \int_{\nu(P,F)} y^p\, dy,$$
where $L(F)$ is the linear vector space parallel to the affine hull of $F$ and $\nu(P,F) \subset S^{n-k-1}$ the set of all outer unit normal vectors to $F\in \mathcal{F}_k(P)$. 

In particular, using the notations from Lemma 4.1 in \citep{HugSchneiderTensorCurvs2013} and identifying $(\R^n)^* \simeq \R^n$:
$$\Phi_{\otimes k,p,1}(K,U) = C_{n,k,p}\, T_K (\mathbf{1}_{(U\times S^{n-1})} \tilde \varphi^{0,p}_k),$$
where $C_{n,k,p}:=(-1)^{n-1}(k-1)!(n-k-1)!\,p!\, s_{n-k+p-1}$ with $s_n:=\vol S^n = \frac{2\pi^{\frac{n+1}{2}}}{\Gamma\left(\frac{n+1}{2}\right)}$.
\end{Prop}

To obtain differential forms with values in an arbitrary irreducible $\SO(n)$-representation $\Gamma_{\lambda}$ from Theorem \ref{thm:curvHarmonicDecomposition}, we need to define two maps.

First, recall that, for any such $\lambda$ with weight $d:=|\lambda|:=\sum_{i=1}^n \lambda_i$, there exists an $\SL(n)$- (hence, also $\SO(n)$-)equivariant projection called the \emph{Young-symmetriser} $\mu_\lambda : (\R^n)^{\otimes d} \to \Gamma_{\lambda}$, where $\Gamma_\lambda$ is the irreducible $\SL(n)$-representation given by $\lambda$. It is best visualised by using the \emph{Young-diagram} associated to $\lambda$, i.e., a left-aligned collection of boxes with $\lambda_i$ boxes in the $i$-th row. The image of $e_{\otimes j_1\ldots j_d}\in (\R^n)^{\otimes d}$ under $\mu_\lambda$ is then represented by the Young-diagram for $\lambda$ with its boxes filled with indexes $j_1,\ldots,j_d$ from top to bottom from left to right. The thus filled diagram is called a \emph{Young-tableau}.

Second, given the canonical projection $\pi_{\tr} : (\R^n)^{\otimes d} \to (\R^n)^{[d]}$ from the $d$-fold tensor product of $\R^n$ to its trace-free subspace, $\bar\Gamma_{[\lambda]}:=\pi_{\tr}(\Gamma_\lambda)$ is an $\SO(n)$-representation. If $n=2m$ is even and $\lambda_m\neq 0$, then $\bar\Gamma_{[\lambda]}$ decomposes into the direct sum of two irreducible $\SO(n)$-representation $\Gamma_{[\lambda]}$ and $\Gamma_{[\bar\lambda]}$, where $\bar\lambda = (\lambda_1,\ldots,\lambda_{m-1},-\lambda_m)$. Otherwise $\bar\Gamma_{[\lambda]}=\Gamma_{[\lambda]}$ is an irreducible $\SO(n)$-representation.

Now, apply $\pi_{\tr} \circ \mu_\lambda$ on the tensor part of the above forms such that the images of $\mu_\lambda$ are given by the following Young-tableaux:
\begin{equation}
	\definecolor{mygray}{rgb}{0.9, 0.9, 0.9}
	\ytableausetup{boxsize=1.55em,aligntableaux=center}
	\FPhi_{[k,p,q]} \sim
	\ytableaushort[*(white)]{ {i_1} {\pi_1} {*(mygray) 1}  {*(mygray) \ldots}   {*(mygray) p},
		\vdots \vdots,
		\vdots \vdots,
		{i_q} {\pi_q}
	},
	\ 
	\FXi_{[k,p,q]} \sim 
	\ytableaushort[*(white)]{ {i_1} {\pi_1} {*(mygray) 1}  {*(mygray) \ldots}   {*(mygray) p},
		\vdots \vdots,
		{i_q} {\pi_q},
		{*(mygray) y} 
	},
	\ 
	\FPsi_{[k,p,q+1]} \sim 
	\ytableaushort[*(white)]{ {i_1} {\pi_1} {*(mygray) 1}  {*(mygray) \ldots}   {*(mygray) p},
		\vdots \vdots,
		{i_q} {\pi_q},
		{*(mygray) y} {*(mygray) y}
	}\hspace{-10pt}
	\label{eq:YoungForms}
\end{equation}
with the integers $j$ in the grey boxes representing the $j$-th copy $y$ in $y^p=y^{\otimes p}$, and symmetrise the tensor part of $\FTheta_{[k,p]}$ as that in $\FPhi_{[k,p,k]}$ except that $\pi_i$ are replaced by $j_i$. We thus obtain the $\bar\Gamma_{[\lambda]}$-valued differential forms:
\begin{equation}
\FT_{[k,p,q]} := \pi_{\tr}\circ \mu_{\lambda} (\FT_{\otimes k,p,q}).
\label{eq:HarmonicForm}
\end{equation}

It is well-known that $\bar\Gamma_{[\lambda]}$ may be embedded into $\AMF{\lambda'}\R^n := \AMF{\lambda'_1}\R^n\otimes \ldots \otimes \AMF{\lambda'_{\lambda_1}}\R^n$, where $\lambda'=(\lambda'_1,\ldots ,\lambda'_{\lambda_1})$ is conjugate to $\lambda$, i.e., where $\lambda'_j$ is the number of boxes in the $j$-th \emph{column} of the Young-diagram of $\lambda$. If $\lambda'_i =n/2$, the operator $*_i: \AMF{\lambda'}\R^n \to \AMF{\lambda'}\R^n$ given by applying the Hodge-$*$-operator on $\AMF{\lambda'_i}\R^n$ restricts to an $\SO(n)$-equivariant map on $\bar\Gamma_{[\lambda]}$ which is not a multiple of the identity. 

\begin{Thm}
\label{thm:OnHarmonicBasis}
Let $\lambda$ be from Theorem \ref{thm:curvHarmonicDecomposition}, $m$ be the largest $j$ with $\lambda_j\neq 0$, $p:=\lambda_1 -2$, and $k':=\min(k,n-k-1)$. Write $\CT_{[k,p,q]}$ for the curvature measure induced by $\FT_{[k,p,q]}$.

$\bullet$ If $m=0$, $\TCurv^{sm, \SO(n)}_{k,\Gamma_{[\lambda]}}$ has the basis $\Phi_{[k,0,0]}$ . 

$\bullet$ If $1\leq m < n/2$, its basis is 
$\begin{cases}
\CXi_{[k,p,m-1]}																& \hbox{if }\lambda_m = 1; 						\\
\CPsi_{[k,p,m]}																	&	\hbox{if }\lambda_m \geq 2 \hbox{ and } m=k'+1;  					\\
\CPhi_{[k,p,m]},\CPsi_{[k,p,m]}\hbox{(, and } \CTheta_{[k,p]}\hbox{)} & \hbox{if }\lambda_m \geq 2 \hbox{ (and } n=2m+1\hbox{)};
\end{cases}
$

$\bullet$ If $m=n/2$, its basis is $\begin{cases}
\CXi_{[k,p,m-1]} \pm i^m *_1\CXi_{[k,p,m-1]}		& \hbox{if }\lambda_m = \mp 1; 						\\
\CPsi_{[k,p,m]} \pm  i^m *_1\CPsi_{[k,p,m]}		&	\hbox{if }\lambda_m = \mp c, c\geq 2.
\end{cases}$

In particular, if $m=n/2$ is odd, $\Gamma_{[\lambda]}$-valued curvature measures cannot be realised as real-valued curvature measures.
\end{Thm}

Although the forms appearing in the above Theorem may seem intimidating at the first glance, they occur naturally when one writes down the isomorphisms in the chain of identities in \pref{eq:HomIsomorphisms} and applies them to the elements of the last space in the chain. The chain itself is the core of the proof of Theorem \ref{thm:curvHarmonicDecomposition} and the elements of the last space are rather straight-forward to construct.

Next, we analyse the behaviour of smooth curvature measures under the globalisation map.

\begin{Thm} 
\label{thm:BernigLemma} 
The kernel of $\glob : \Curv^{sm}_k \to \Valsm_k$ is spanned by:
	\begin{eqnarray}
	\CXi_{[k,p,q]}, \CTheta_{[k,p]} & & \quad \hbox{for all } p,q, \label{eq:GlobXi} \\
	\CPsi_{[k,p,k+1]} 		& & \quad \hbox{for all $p$ if $0\leq k\leq \frac{n-1}{2}$} \label{eq:GlobPsiMargin} \\
	q (n-k+1)\, \CPsi_{[k,p,q]} + (k-q+1)(q p + 1)\, \CPhi_{[k,p,q]}  & &  \quad  \hbox{for all $p$ and $1 \leq q\leq k'$}. \label{eq:GlobPhiPsi}
	\end{eqnarray}
\end{Thm}

This yields in combination with Proposition \ref{thm:smoothness} the following result. 

\begin{Prop}
\label{thm:TValBasis}
Let $\lambda$ be from the harmonic decomposition of $\Val_k$ such that all $\lambda_j \geq 0$. Writing $\tau_{[k,p,q]} := \glob T_{[k,p,q]}$, the space $\TVal^{\SO(n)}_{k,\bar\Gamma_{[\lambda]}}$ is spanned by $\VPhi_{[k,0,0]}$ if $m = 0$, $\VPsi_{[k,p,m]}$ if $\lambda_m \geq 2$ and $m<n/2$, and $\VPsi_{[k,p,m]},*_1\VPsi_{[k,p,m]}$ otherwise.

In particular, the coefficients of $\VPhi_{[k,0,0]}$, $\VPsi_{[k,p,q]}$, $1\leq q \leq \min\{k,n-k\}, p\geq 0$ -- and those of $*_1\VPsi_{[k,p,k]}$ if $n=2k$ -- form a Schauder-basis of $\Val_k$.
\end{Prop}

\section{Representation Theory}
\label{sec:RepTheory}

Let $V=\C^n$ with $n\geq 3$ and assume that all representations are finite-dimensional -- unless otherwise stated -- in this section.

Given a Young-diagram $\lambda$, define two subgroups of the permutation group $S_d$:
\begin{eqnarray*}
P & = & \{\pi\in S_d:\ \pi \hbox{ preserves each row of $\lambda$}\}, \\
Q & = & \{\pi\in S_d:\ \pi \hbox{ preserves each column of $\lambda$}\}.
\end{eqnarray*}
Defining the \emph{group algebra} $\C G$ to be a vector space spanned by vectors $e_g$ for each $g\in G$, such that $e_g\cdot e_h = e_{gh}$, we set:
\begin{equation}
a_\lambda = \sum_{\pi\in P} e_\pi \in \C S_d, \quad\quad b_\lambda = \sum_{\pi\in Q} \sgn \pi\cdot e_\pi \in \C S_d, \quad \hbox{ and } \quad c_\lambda = a_\lambda \cdot b_\lambda \in \C S_d.
\label{eq:YoungSymmetrisers}
\end{equation}
It turns out that $c_\lambda \cdot c_\lambda = n_\lambda c_\lambda$ for some positive integer $n_\lambda$ and $\SS_\lambda V := V^{\otimes d}\cdot c_\lambda$ is an irreducible $S_d$-representation. Furthermore, the right action of $S_d$ on $V^{\otimes d}$ given by permuting factors $(v_1\otimes\dots\otimes v_d)\cdot \sigma = v_{\sigma(1)}\otimes\cdots \otimes v_{\sigma(d)}$
commutes with the standard left action of $\SL(n,\C)$. Hence, $\SS_\lambda V$ is also an irreducible $\SL(n,\C)$-module. The map $\mu_\lambda(v):= v\cdot c_\lambda$ is the \emph{Young-symmetriser} mentioned in the introduction. 

\begin{Prop}
Any irreducible complex $\SL(n,\C)$-module is isomorphic to the $\SL(n,\C)$-module $\SS_\lambda V$ for some $\lambda=(\lambda_1\geq \ldots \geq \lambda_n\geq 0)$. The isomorphy class of $\SL(n,\C)$-representations which contains $\SS_\lambda V$ is denoted by $\Gamma_\lambda$.
\label{thm:SLrepClassification}
\end{Prop}

See \citep[Chapter 6.1, Proposition 15.15]{FultonHarris} for more details.

The $\SL(n,\C)$-modules $\Gamma_\lambda$ are also uniquely determined up to isomorphism by certain Bianchi-type identities \citep[\S 8]{Fulton}, \citep[\S I.5, (5.12)]{Macdonald}. Define a \emph{(Young)-tableau} $T$ on $\lambda$ as a numbering of the boxes by the integers $1,\ldots, |\lambda|=:d$ and let $T(i,j)$ be the number in the $i$-th box of the $j$-th column. A \emph{semi-standard tableau} is a Young-tableau such that the entries are non-decreasing in each row and strictly increasing in each column.

\begin{Thm}[Bianchi-type identities]
\label{thm:BianchiRelationsSAV}
Let $e_1,\ldots, e_n$  be an orthonormal base of $V$ and write $e_T:= \prod_{j=1}^{\lambda_1} e_{T(1,j)}\otimes \cdots \otimes e_{T(\lambda'_j,j)}\in V^{\otimes |\lambda|}$ for any Young-tableau $T$ of $\lambda$. Then for any \emph{semi-standard} tableau $T$, one has:
$$\textstyle\mu_{\lambda}\left(e_T - \sum_S e_S\right) = 0,$$	
where the sum is over all $S$ obtained from $T$ by exchanging the top $k$ elements of one column with any $k$ elements of the preceding column, maintaining the vertical orders of each set exchanged. There is one such relation for each numbering $T$, each choice of adjacent columns, and each $k$ at most equal to the length of the shorter column. 

The elements $\mu_{\lambda}(e_T)$ for semi-standard Young-tableaux $T$ generate $\SS_\lambda V$ as a vector space.
\end{Thm}

$\SS_\lambda V$ may be used to construct irreducible $\SO(n,\C)$- and $\OLie(n,\C)$-modules. As there exists a symmetric bilinear form $Q$ on $V$ preserved by $\OLie(n,\C)$, the contraction maps for $p<q$:
\begin{eqnarray}
\tr_{p,q} :\qquad\quad V^{\otimes d}\ \quad& \to & V^{\otimes d-2} \label{eq:traceSOn} \\
				v_1\otimes\cdots\otimes v_d & \mapsto & Q(v_p, v_q)\,v_1\otimes\cdots\otimes \hat v_{p}\otimes\cdots\otimes \hat v_{q} \otimes \cdots \otimes v_d \nonumber
\end{eqnarray}
are $\OLie(n)$-equivariant. The intersection of all kernels of such contractions is closed under the action of $S_d$, hence, the intersection $V^{[d]}$ of these kernels is an $S_d$-submodule of $V^{\otimes d}$. Set
$
\SS_{[\lambda]}V:= V^{[d]} \cap \SS_\lambda V.
$

\begin{Thm}
The $\OLie(n,\C)$-module $\SS_{[\lambda]} V$ is trivial if $\lambda_{\floor{n/2}+1}> 0$ or $\lambda'_1+\lambda'_2 > n$ and irreducible otherwise. Furthermore:
\begin{itemize}
	\item If $n=2k+1$ and $\lambda = (\lambda_1\geq\lambda_2\geq\ldots\lambda_k\geq 0)$ or  $n=2k$ and $\lambda = (\lambda_1\geq\lambda_2\geq\ldots\lambda_{k-1}\geq \lambda_k=0)$,
then $\SS_{[\lambda]}V$ is an irreducible $\SO(n,\C)$-representation.
	\item If $n = 2k$ and $\lambda = (\lambda_1\geq\lambda_2\geq\ldots\lambda_k > 0)$, then $\SS_{[\lambda]} V$ is a direct sum of two irreducible $\SO(n,\C)$-modules that are dual to each other.
\end{itemize}
\label{thm:SOnRepClassification}
\end{Thm}

We write $\bar\Gamma_{[\lambda]}$ for the isomorphy class of irreducible $\OLie(n, \C)$-representations containing $\SS_{[\lambda]} V$ and $\Gamma_{[\lambda]}$ for the isomorphy class of irreducible $\SO(n,\C)$-representations corresponding to the tuple $\lambda$. One may show that $\Gamma_{[\lambda_1,\ldots,\lambda_k]}^* = \Gamma_{[\lambda_1,\ldots,\lambda_{k-1},-\lambda_k]}$ and the theorem may be re-stated as: 
\begin{equation*}
\bar\Gamma_{[\lambda]} =  
\begin{cases} 
\Gamma_{[\lambda]} \oplus \Gamma^*_{[\lambda]} & \quad \hbox{if $n=2k$ is even and $\lambda_{k}\neq 0$,} \\
\Gamma_{[\lambda]} & \quad \hbox{otherwise}.
\end{cases}
\end{equation*}

\begin{Def}
Let $V$ be a representation of a Lie-group $G$. The \emph{character} $\Char{V}$ of $V$ is a complex-valued function on $G$ defined by $\chi_V(g)=\mathrm{tr}(g|_V)$.
\end{Def}
The most notable facts about characters is their ability to uniquely determine $G$-modules up to isomorphism for any compact or linear reductive Lie-group $G$ as well as their explicit forms for a large number of representations. For example, the character of the irreducible $\SL(n,\C)$-module $\AMF{k}V$ is given by the elementary symmetric polynomial $E_k$ of the eigenvalues $x_1,\ldots, x_n$ of $g\in \SL(n,\C)$:
\begin{equation*}
\Char{\AMF{k}V} (g) = E_k(x_1,\ldots,x_n) \,= \sum_{i_1<\cdots<i_k=1}^n \!x_{i_1}\cdot \ldots \cdot x_{i_k}.
\end{equation*}
More generally, one has the following result.
\begin{Prop}[Giambelli-formula for $\SL(n,\C)$] \label{thm:detFormula}
Let $\lambda$ be a tuple $(\lambda_1\geq\ldots \geq \lambda_n\geq 0)$ and $\mu = (\mu_1,\ldots,\mu_\ell) = \bar\lambda$ its conjugate partition. Then:
\begin{equation*}
\Char{\Gamma_\lambda} = \det ( E_{\mu_i + j - i} ) = 
\det \begin{pmatrix} 
E_{\mu_1} & E_{\mu_1 + 1} & \cdots & E_{\mu_1 + \ell - 1} \\
E_{\mu_2 -1 } & E_{\mu_2} & \cdots & E_{\mu_2 + \ell - 2} \\
\vdots & \vdots & \ddots & \vdots \\
E_{\mu_\ell-l+1} & E_{\mu_\ell-l} & \cdots & E_{\mu_\ell}
\end{pmatrix}.
\end{equation*}
\end{Prop}

A similar formula may be found for characters of $\SO(n,\C)$-representations except that the character of $\AMF{k}V$ as an $\SO(n)$-representation is given by $E_k=E_k(x_1,\ldots,x_m,x_1^{-1},\ldots, x_m^{-1})$ for $n=2m$ and $E_k = E_k(x_1,\ldots,x_m,x_1^{-1},\ldots, x_m^{-1},1)$ for $n=2m+1$. Then $E_{m+k}=E_{m-k}$ resp. $E_{m+k} = E_{m+1-k}$ due to the isomorphisms $\AMF{m+k}V \simeq \AMF{m-k}V$ resp. $\AMF{m+k}V \simeq \AMF{m-k+1}V$ for even resp. odd $n$. 

\begin{Prop}[Giambelli-formula for $\SO(n,\C)$]\label{thm:detFormulaSOn}
Let $\lambda$ be a tuple of integers $(\lambda_1\geq\ldots \geq \lambda_n\geq 0)$ and $\mu = (\mu_1,\ldots,\mu_\ell) = \bar\lambda$ its conjugate partition. Then the character $\Char{\bar\Gamma_{[\lambda]}}$ is given by the determinant of the $\ell\times\ell$-matrix with $i$-th row
$$ (E_{\mu_i - i + 1} \quad E_{\mu_i - i + 2} + E_{\mu_i - i} \quad E_{\mu_i - i + 3} + E_{\mu_i - i -1} \quad \cdots \quad E_{\mu_i - i + \ell} + E_{\mu_i - i - \ell + 2}).$$
\end{Prop}

Given a representation $V$ of a Lie-group $G$, any closed Lie-subgroup $H\subset G$ inherits from $G$ the action on $V$ so that $V$ may also be regarded as an $H$-module which we denote by $\Res^G_H V$. 
Such restrictions may often be written in closed terms. 
\begin{Thm}[$\SO(n,\C)$-branching]
Let $\lambda$ be a tuple of integers satisfying conditions from Theorem \ref{thm:SOnRepClassification}. Then
$$
\Res_{\SO(n-1,\C)}^{\SO(n,\C)}\Gamma^{\SO(n,\C)}_{[\lambda]}=\bigoplus_{\mu} \Gamma^{\SO(n-1,\C)}_{[\mu]},
$$
where $\mu$ runs over all partitions $\mu=(\mu_1,\ldots,\mu_{k})$, $k=\floor*{(n-1)/2}$, such that
$$
\begin{cases}
\lambda_1\geq\mu_1\geq\lambda_2\geq\mu_2\geq\ldots\geq\mu_{k-1}\geq\lambda_{\floor*{n/2}}\geq|\mu_k| & \hbox{ for odd } n, \\
\lambda_1\geq\mu_1\geq\lambda_2\geq\mu_2\geq\ldots\geq\mu_{k} \geq|\lambda_{\floor*{n/2}}| & \hbox{ for even } n. 
\end{cases}
$$
\label{thm:BranchingSOn}
\end{Thm}

There is also a canonical way to ``extend" a representation $W$ of $H$ to a representation of $G$. Consider the space $C^\infty(G,W)$ of all smooth functions from $G$ to $W$. The $G$-invariant subspace:
\begin{equation}
\Ind^G_H W :=\{f\in C^\infty(G,W)\,|\, f(gh) = h^{-1} f(g),\quad \forall h\in H,\, \forall g\in G\}.
\label{eq:analyticInd}
\end{equation}
is called the induced representation of $G$ from $H$. 

Note that $\Ind^G_H W$ is, in general, not finite-dimensional. Nevertheless, the formulae for $\Res\left(\Ind W\right)$ and $\Ind\left(\Res W\right)$ are known and can be found in \citep{Serre}. Although both constructions are generally not equal to $W$, the well-known Frobenius' Theorem shows that $\Ind$ and $\Res$ are, in some sense, adjoint to each other. 

\begin{Thm}[Frobenius' Reciprocity Theorem]
\label{thm:Frobenius}
Let $G$ be a compact Lie-group and $H\subset G$ a closed Lie-subgroup. Given a representation $U$ of $G$ and a representation $W$ of $H$, there is a canonical vector space isomorphism 
$$\Hom_G(U, \Ind W) \simeq \Hom_H(\Res U, W).$$
\end{Thm}

We can now prove the following result which is a refinement of Corollary 3.4 in \citep{AleskerBernigSchuster}.
\begin{Lemma}
\label{thm:LambdaDecomposition}
Let $i,j\in\N$ such that $0\leq i, j\leq n$ and set 
$$i':=\max(\min(i,n-i),\min(j,n-j)), \qquad j':=\min(\min(i,n-i),\min(j,n-j)).$$
Then the following $\SL(n,\C)$-representations are isomorphic:
\begin{equation}
	\AMF{i,j} V \simeq \left(\Gamma_{(2[j'], 1[i'-j'])}\right) \oplus \AMF{i'+1, j'-1} V \simeq \bigoplus_{k=0}^{j'} \bar\Gamma_{(2[j'-k], 1[2k+i'-j'])}
\label{eq:LambdaDecomposition1}
\end{equation} 
The above isomorphisms may be interpreted as isomorphisms of $\SO(n,\C)$-representations by the following identity of $\SO(n)$-representations $\Res \Gamma_{(2[k],1[l])} = \bigoplus_{m=0}^k \Gamma_{[2[m],1[l]]}$ for any integers $k, l$. 
\end{Lemma}
\begin{proof} 
Since $\AMF{i}V\simeq \AMF{n-i}V$ and $\AMF{i}V\otimes\AMF{j}V \simeq \AMF{j}V\otimes\AMF{i}V$, we may assume w.l.o.g. $i=i' \leq n/2$ and $j=j'\leq n/2$. If $\lambda=(\lambda_1,\ldots,\lambda_{m})$ is a non-negative tuple, as specified in the middle term of the above identity, then the conjugate $ \mu :=\lambda' =(i, j)$. By Proposition \ref{thm:detFormula}:
$$\Char{\Gamma_\lambda} = \det
\begin{pmatrix}
	E_{i}  & E_{i+1} \\
	E_{j-1} & E_{j}
\end{pmatrix} = E_i E_j - E_{i+1} E_{j-1},
$$
which shows the left isomorphism in \pref{eq:LambdaDecomposition1}. Applying it recursively until $j'=0$ yields the right isomorphism. Apply Proposition \ref{thm:detFormulaSOn} on $\bar \Gamma_{[\lambda]}$ for $\lambda = (2[m],1[l])$ with conjugate $\mu=(l+m,m)$:
$$ \Char{\bar\Gamma_{[\lambda]}} = \det
\begin{pmatrix}
	E_{m + l} & E_{m+l+1} + E_{m+l-1} \\
	E_{m - 1} & E_m + E_{m -2}
\end{pmatrix}.
$$
The last identity is now obtained by summing over all $m$:
\begin{eqnarray*}
\sum_{m=0}^k \Char{\bar\Gamma_{[2[m],1[l]]}} & = & \sum_{m = 0}^k \left( E_{m+l}(E_{m} + E_{m - 2}) - E_{l -1}(E_{m+l+1} + E_{m+l-1})\right) \\
 & = & E_{k+l} E_k - E_{k+l+1}E_{k-1} = \Char{\bar\Gamma{(2[k],1[l])}}.
\end{eqnarray*}
\vspace{-.9cm}

\end{proof}

\begin{Remark}
The complexification of $\so(n,\R)$ is $\so(n,\C)$ and that of $\sl(n,\R)$ is $\sl(n,\C)$ which are both complex simple Lie-algebras. By \citep[Chapter 5.1]{Knapp}, \citep[Chapter 26.1]{FultonHarris}, if $G$ is a real Lie-group with a simple real Lie-algebra $\mg_0$ such that its complexification $\mg:=\mg_0\otimes \C$ is a simple complex Lie-algebra, then there is one-to-one correspondence between the complex representations of $G$ and its complexified counterpart with the Lie-algebra $\mg$. Thus, one obtains a one-to-one correspondence between the complex representations of $\SO(n):=\SO(n,\R)$ resp. $\SL(n):=\SL(n,\R)$ and those of $\SO(n,\C)$ resp. $\SL(n,\C)$. 
\end{Remark}
\begin{Remark} \label{rem:realModules}
The $\SO(n)$-module $\Gamma_{[\lambda]}$ on a complex vector space is called \emph{of real type} (or just real) if it may be realised as a complexification $\Gamma_{[\lambda,\R]}\otimes \C $ of an irreducible $\SO(n)$-module with the same tuple $\lambda$ on real vector space. By \citep[Proposition 26.27]{FultonHarris}, the $\SO(n)$-module $\Gamma_{[\lambda]}$ is not of real type if and only if $n=2k$ for odd $k$ and $\lambda_k\neq 0$. In contrast, irreducible $\OLie(n)$-modules $\bar\Gamma_{[\lambda_1,\ldots,|\lambda_k|]}$ are always of real type.
\end{Remark}

\section{Valuation Theory and Contact Geometry}
\label{sec:valuations}
From now on, we assume that $V=\R^n$ with the basis $e_1,\ldots,e_n$ and write $\SL(n)=\SL(n,\R)$ and $\SO(n) = \SO(n,\R)$. 

The \emph{normal cycle} of a convex body $K\in \KK(V)$ is an $(n-1)$-dimensional Lipschitz manifold:
$$
\nc(K) := \{(x,y)\in SV\, |\, \skprod{x-x',y}\geq 0,\,\forall x'\in K\}.
$$
\begin{Def}
\label{def:smoothVal}
A translation-invariant functional $\VPhi : \KK(V)\to A$ is called a \emph{smooth valuation} if, for all $K\in\KK(\R^n)$,
\begin{equation*}
\VPhi(K) = \integ(\beta, \omega)(K) := \int_K \beta + \int_{\nc(K)} \omega,
\end{equation*}
where $\beta\in\SDF{n}{\R^n}{\R^n}\otimes \Gamma$ is a translation-invariant $\Gamma$-valued form on $\R^n$ and $\omega\in\SDF{n-1}{S\R^n}{\R^n}\otimes \Gamma$ is a translation-invariant form on $S\R^n$. Likewise, a translation-invariant functional $\CPhi : \KK(V)\times \BB(V) \to \Gamma $ is called a \emph{smooth curvature measure} if, for all $K\in\KK(\R^n)$ and all $U\in\BB(\R^n)$, 
\begin{equation*}
\CPhi(K,U) = \Integ(\beta, \omega)(K,U) := \int_{K\cap U} \beta + \int_{\nc(K)\cap \pi^{-1}(U)} \omega,
\end{equation*}
where $\pi:S\R^n \to \R^n$ is the projection on the first factor. The operators $\integ$ and $\Integ$ which assign to a given pair of translation-invariant forms a corresponding smooth valuation resp. curvature measure are called the \emph{integration operators}.  
\end{Def}

Both integration operators have non-trivial kernels best described in contact-geometric terms. Let $(W,\omega)$ be a symplectic vector space of real dimension $2n$. Recall that the operator 
\begin{eqnarray*}
\Lef : \AMF{*}(W^*) & \rightarrow & \AMF{* +2}(W^*) \nonumber \\
          \tau & \mapsto     & \tau \wedge \omega
\end{eqnarray*}
is called the \emph{Lefschetz operator}. Fixing an Euclidean scalar product $\skprod{\cdot,\cdot}$ on $W$, the operator $\DLef$ of degree $(-2)$ uniquely determined by
\begin{equation*}
\skprod{\DLef\tau, \beta} = \skprod{\tau,\Lef\beta},\quad \forall\beta,\tau\in\AMF{*}(V^*)
\end{equation*}
is called the \emph{dual Lefschetz operator}.
\begin{Def}
A $k$-linear form $\alpha\in\AMF{k}(W^*)$ is called \emph{primitive} if $\DLef\alpha=0$. The subspace of all primitive elements in $\AMF{k}(W^*)$ is denoted by $\AMFP{k}(W^*)\subset\AMF{k}(W^*)$. The operator $\DLef$ and, hence, the notion of primitivity may be extended to \emph{symplectic manifolds} in a pointwise manner.
\end{Def}

To define a contact manifold, recall that a \emph{contact element} on a manifold $M$ is a point $p\in M$, called the contact point, together with a tangent hyperplane at $p$, $Q_p\subset T_p M$, i.e. a co-dimension 1 subspace of $T_p M$. A hyperplane $Q_p\subset T_p M$ is completely determined by a linear form $\alpha_p\in T^*_p M \setminus\{0\}$ that is unique up to some non-zero scalar. Indeed, if $(p,Q_p)$ is a contact element, then $Q_p=\ker\alpha_p$. On the other hand, $\ker \alpha_p = \ker \alpha'_p$ if and only if $\alpha_p=\lambda \alpha'_p$. Now, let $Q$ be a smooth field of contact hyperplanes on $M$ defined by $Q(p) := Q_p$. Then $Q=\ker\alpha$ for an open subset $U \subset M$ and some $1$-form $\alpha$ called a locally defining $1$-form for $Q$. This form is again unique up to a smooth nowhere vanishing function $f\in C^\infty(U)$.

A \emph{contact structure} on $M$ is a smooth field of tangent hyperplanes $Q\subset TM$ such that, for any locally defining $1$-form $\alpha$, $\local{d\alpha}{Q}$ is non-degenerate, i.e. symplectic. The pair $(M,Q)$ is called a \emph{contact manifold} and $\alpha$ is called a \emph{local contact form}. The restriction $\local{d\alpha_p}{Q_p}$ is symplectic on $Q_p$,  which implies immediately that $\dim Q_p= 2n$ is even and $\local{d\alpha_p^n}{Q_p}\neq 0$ is a volume form on $Q_p$. Since $T_p M = \ker\alpha_p\oplus\ker d\alpha_p$, one has $\dim T_p M = 2n+1$ is odd. In fact, $Q$ is a contact structure if and only if $\alpha\wedge d\alpha^n\neq 0$ for every locally defining $1$-form $\alpha$. In particular, $\alpha$ is a global contact form if and only if $\alpha\wedge d\alpha^n$ is a volume form on $M$

If there is a globally defined form $\alpha$, one can obtain a unique vector field $T$ called the \emph{Reeb vector field} on $M$ such that the contraction $\iota_T(d\alpha)=0$ and $\iota_T(\alpha)=1$. Indeed, $\iota_T(d\alpha)=0$ implies that $T\in\ker d\alpha$, which is one-dimensional, and $\iota_T{\alpha}=1$ just normalises $T$. 

We may now refine the description of differential forms on $S\R^n$ which turns out to be a contact manifold with the contact form $\alpha$ defined pointwise at $p=(x,y)\in S\R^n$ as follows:
$$\local{\alpha}{(x,y)}(w):=\skprod{y,d\pi(w)} = \sum_{i=1}^n y_i\, \dxi{i}(w),$$
where $\pi:S\R^n\rightarrow \R^n$ is the projection. The Reeb vector field $T$ is given by $\local{T}{(x,y)}=\sum_{i=1}^n y_i \pdxi{}{i}$.

\begin{Def}
A form $\omega\in\Omega^*(S\R^n)$ is called \emph{horizontal} if $\iota_T{\omega}=0$. A form $\omega$ that can be written as $\tau\wedge\alpha$ is called \emph{vertical}. The algebras of horizontal or vertical forms on $S\R^n$ are denoted by $\SDFH{*}{S\R^n}{}$ and $\SDFV{*}{S\R^n}{}$, respectively.
\end{Def}

A smooth translation-invariant form $\omega$ on $S\R^n$ is said to be of bi-degree $(i,j)$ if $\omega$ can be written as $\sum_a\tau_a \otimes\phi_a$ with $\tau_a\in\SDF{i}{\R^n}{\R^n}$ and $\phi_a\in\SDF{j}{S^{n-1}}{}$. Clearly, $\omega\in\SDF{i+j}{S\R^n}{\R^n}$ and
\begin{equation*}
\SDF{k}{S\R^n}{\R^n} = \bigoplus_{i+j=k} \SDF{i}{\R^n}{\R^n} \otimes \SDF{j}{S^{n-1})}{}.
\end{equation*}
To simplify the notation, we write $\Omega^{i,j}$ for the space $\SDF{i,j}{S\R^n}{\R^n}$ of translation-invariant differential forms of bi-degree $(i,j)$ on $S\R^n$ and $\Omega_p^{i,j}$ for the space of \emph{primitive} translation-invariant forms. As $\alpha\in\Omega_v^{1,0}$ and $L$ is of bi-degree $(1,1)$ in this notation, we have:
\begin{equation}
\Omega_p^{i,j} = \Omega_h^{i,j} / L \Omega_h^{i-1,j-1},
\label{eq:primitiveFormsAsQuotient}
\end{equation}
whenever $i+j \leq n$. Furthermore, the Hodge-$*$-operator on $S\R^n$ induces two finer operators on $\Omega^{*}$: $*_1 : \Omega^{i,j} \rightarrow \Omega^{n-i,j}$ and $*_2:\Omega^{i,j} \rightarrow \Omega^{i,n-j-1}$ given by applying the Hodge-$*$-operator on the $\SDF{i}{\R^n}{\R^n}$- resp. $\SDF{j}{S^{n-1}}{}$-part of a differential form. Since, for any vertical translation-invariant form $\omega$, both $*\omega$ and $*_1\omega$ are translation-invariant and horizontal, and vice versa, both operators yield isomorphisms $ *_1 : \Omega_h^{i,j} \to \Omega_h^{n-1-i,j}$ and $*_2:\Omega_h^{i,j} \rightarrow \Omega_h^{i,n-j-1}$. 

To reduce a vertical form $\tau\wedge\alpha$ to a horizontal form, we use a contraction with the Reeb vector field $\iota_T$. Indeed, $\iota_T(\tau\wedge\alpha) = (\iota_T\tau)\wedge\alpha + \tau\wedge(\iota_T\alpha) = \tau$ for any horizontal form $\tau$. Hence, we may write for $\omega\in\Omega^{i,j}$ (recall that $\AMF{i,j}\R^n = \AMF{i}\R^n\otimes \AMF{j}\R^n$):
$$\omega|_{(x,y)} \in \left(\AMF{i,j} T_y^{*}S^{n-1}\right) \oplus \left(\AMF{i-1,j} T_y^{*}S^{n-1}\otimes \R \alpha|_{(x,y)}\right). $$
In particular, if $\omega\in\Omega_h^{i,j}$, then $\omega|_{(x,y)} \in \AMF{i,j} T_y^{*}S^{n-1}.$ We will write in the following $\omega|_y$ instead of $\omega|_{(x,y)}$, whenever $\omega\in\Omega_h^{i,j}$ and $(x,y)\in S\R^n$. Observing that the stabiliser of $\SO(n)$ at any fixed point $y\in S^{n-1}$ is $\SO(n-1)$ and writing $W_y := T_yS^{n-1}$, one has the following result.

\begin{Lemma}[\citep{AleskerBernigSchuster}]
\label{thm:horizontalInduced} 
For all $i,j\in\N$, one has $\Omega_h^{i,j} \simeq \Ind_{\SO(n-1)}^{\SO(n)}(\AMF{i,j} W_y^*).$
\end{Lemma}
\begin{Cor}
\label{thm:primitiveDecomposition}
If $i+j\leq n-1$ and $\max(i,j)\geq (n-1)/2$, then there is an isomorphism of $\SO(n)$-representations 
\begin{equation}
\label{eq:primitiveDecomposition}
\Omega_p^{i,j}{} \oplus \Ind_{\SO(n-1)}^{\SO(n)}(\AMF{i-1,j-1} W_y^*) = \Ind_{\SO(n-1)}^{\SO(n)}(\AMF{i,j} W_y^*),
\end{equation}
hence, $\Omega_p^{i,j} = \Ind_{\SO(n-1)}^{\SO(n)}\AMFP{i,j} W^*_y$, where $\AMFP{i,j} W^*_y := \bigoplus_{l=0}^j \bar \Gamma_{[2[l], 1[n-1-(i+j)]]}$. 
\end{Cor}
\begin{proof}
Let w.l.o.g. $j \geq (n-1)/2$. Then $i\leq (n-1)/2$ and Lemma \ref{thm:LambdaDecomposition} yields:
\begin{eqnarray*}
\AMF{i,j}W^*_y & = & *_2(\AMF{i,n-j-1}W^*_y) 
		 \simeq   *_2(\AMFP{i,n-j-1} W^*_y) \oplus *_2(\AMF{\min\{i,n-j-1\}-1, \max\{i,n-j-1\}+1} W^*_y)
\end{eqnarray*}
As $W^*_y\oplus W^*_y$ is a symplectic space with the symplectic form $d\alpha$ and 
$$ *_2(\AMF{\min\{i,n-j-1\}-1,\max\{i,n-j-1\}+1} W^*_y) \subset \AMF{i+j-2} (W_y^*\oplus W_y^*),$$
the Lefschetz decomposition implies that
$$\AMF{i,j}W^*_y = *_2(\AMFP{i,n-j-1} W^*_y) \oplus d\alpha \wedge *_2(\AMF{\min\{i,n-j-1\}-1,\max\{i,n-j-1\}+1} W^*_y).$$
The claim follows now immediately from \pref{eq:primitiveFormsAsQuotient} and the above Lemma. Note that the condition $\max(i,j)\geq (n-1)/2$ is essential for the claim's validity.
\end{proof}
 
\begin{Thm}
The $\SO(n)$-representations $\Curv^{sm}_k$ and $\Omega_p^{k,n-1-k}$ are isomorphic and one has:
\label{thm:CurvIsomorphicToPrimitive}
\begin{equation*}
\Curv^{sm} = \bigoplus_{k=0}^n \Curv^{sm}_k.
\end{equation*}
\end{Thm}
\begin{proof}
We know from \citep{BernigLefschetz} that $\ker \integ$ is generated by vertical and exact forms and it is obvious that $\ker \Integ \subset \ker \integ$. Vertical forms are precisely those which vanish pointwise on normal cycles, hence, they lie in $\ker \Integ$. Let $\omega = d\tau$ be an exact horizontal $(n-1)$-form. Then, for $K\in\KK(\R^n)$ and $U\in\BB(\R^n)$:
\begin{equation}
\int_{\nc(K)\cap \pi^{-1}(U)} d\tau = \int_{\partial(\nc(K)\cap \pi^{-1}(U))} \tau.
\label{eq:intExactVanishes}
\end{equation}

Since $\partial(\nc(K)\cap \pi^{-1}(U)) \subset \nc(K)$, the integral vanishes for any $K$ and $U$ if and only if $\tau$ vanishes on $\nc(K)$ pointwise, i.e., if $\omega = d(\alpha\wedge\phi) = d\alpha\wedge \phi - \alpha \wedge d\phi$. The second term is 0 due to horizontality of $\omega$, hence, $\omega$ is a multiple of $d\alpha$ and the first claim follows. The decomposition of $\Curv^{sm}$ follows immediately from the bi-grading on $\Omega_p^{*}$.

\end{proof}
\section{Proofs of the Main Results}
\label{sec:coreResults}

\subsection{Decomposition and Basis}
\label{sec:HarmonicDecomposition}
A tuple $\lambda$ is said to be of type $[q;p;r]$ if its conjugate is $(q+r,q,1\ldots,1)$ and $q+r$ or $q$ are ignored if they are 0. The $\SL(n)$- and $\SO(n)$-representations associated to such tuples are also called of type $[q;p;r]$. In particular, the representation of type $[0;0;0]$ is trivial and that of types $[0;0;1]$ or $[0;1;0]$ is the standard representation. Theorem \ref{thm:curvHarmonicDecomposition} claims that only $\SO(n)$-representations $\Gamma_{[\lambda]}$ of type $[q;p;r]$ and their duals occur in $\Curv^{sm}_k$. These $\SO(n)$-representations will be denoted by $\Gamma^{q,p}_r$.

\begin{proof}[Proof of Theorem \ref{thm:curvHarmonicDecomposition}] We write $n':=n-1$ for brevity and assume w.l.o.g $k\leq n'/2$. To distinguish between $\SO(n)$- and $\SO(n')$-representations, we denote the former by $\Gamma_{[\lambda]}$ and the latter -- by $\Upsilon_{[\lambda]}$. The operators $\Res^{\SO(n)}_{\SO(n')}, \Ind^{\SO(n)}_{\SO(n')}$ will be shortened to $\Res$ and $\Ind$, respectively. 

Let $\Gamma_{[\lambda]}$ be an arbitrary irreducible $\SO(n)$-representations. By Schur's Lemma, the total multiplicity of $\Gamma_{[\lambda]}$ in $\Curv^{sm}_k$ is the dimension of $\Hom_{\SO(n)} (\Curv^{sm}_k, \Gamma^*_{[\lambda]})$. As $\Hom_G(V,W)^G \simeq (V^*\otimes W)^G$, one has:
\begin{equation*}
(\Curv^{sm}_k \otimes \Gamma_{[\lambda]})^{\SO(n)} \overset{Thm\ \ref{thm:CurvIsomorphicToPrimitive}}{=} (\Omega_p^{k,n'-k} \otimes \Gamma_{[\lambda]})^{\SO(n)}  
	\overset{Cor\ \ref{thm:primitiveDecomposition}}{=} \left(\Ind \Lambda_p^{k,n'-k} W^*_y \otimes \Gamma_{[\lambda]}\right)^{\SO(n)} 	
\end{equation*}																						
\begin{equation}
	\mkern-48mu \overset{Thm\ \ref{thm:Frobenius}}{=} \bigoplus_{q=0}^k \Hom_{\SO(n')}\left(\bar \Upsilon^{q,0}_0, \Res \Gamma_{[\lambda]}\right) 
	= \bigoplus_{q=0}^k \bigoplus_\mu \Hom_{\SO(n')}\left(\bar \Upsilon^{q,0}_0, \Upsilon_{[\mu]}\right),
\label{eq:HomIsomorphisms}
\end{equation}
where the sum over $\mu$ is as per Theorem \ref{thm:BranchingSOn}. Note that we have dropped the duality in the third equality, since $\Res \Gamma_{[\lambda]} \simeq \Res (\Gamma_{[\lambda]})^* \simeq (\Res \Gamma_{[\lambda]})^*$ and, hence, the multiplicity of $\Gamma_{[\lambda]}$ and $(\Gamma_{[\lambda]})^*$ in $\Curv_k^{sm}$ is the same. By Schur's Lemma, $\Hom_{\SO(n')}(\bar \Upsilon^{q,0}_0, \Upsilon_{[\mu]})$ is not trivial if and only if $\mu=[q;0;0]$. Hence, the multiplicity of $\Gamma_{[\lambda]}$ in $\Curv^{sm}_k$ is equal to the number of modules of type $[q;0;0]$ in $\Res\Gamma_{[\lambda]}$. We now study the classes of $\Gamma_{[\lambda]}$ on a case-by-case basis:
\begin{itemize}
	\item $\Gamma^{q,p}_1$ contains exactly one $\SO(n')$-module $\Upsilon^{q,0}_0$ if and only if $0\leq q\leq k$.
	\item $\Gamma^{q,p}_0$ contains modules $\bar\Upsilon^{q,0}_0, \bar\Upsilon^{q-1,0}_0$ if $1\leq q\leq k$, $\bar\Upsilon^{k,0}_{0}$ if $q=k+1$, and $\bar\Upsilon^{0,0}_0$ if $q=p=0$. Note that $\bar\Upsilon^{q,0}_0$ is a sum of two irreducible modules if and only if $q=k=n'/2$, i.e., when $n=2k+1$, otherwise it is irreducible.
	\item The same applies for the above modules' duals. The only non-self-dual modules with non-zero multiplicities in $\Curv_k$ are $(\Gamma^{k,p}_{0})^*$ and $(\Gamma^{k-1,p}_1)^*$ if $n=2k$.
\end{itemize}
Irreducible $\SO(n)$-modules not mentioned in the above list do not contain $\SO(n-1)$-modules of type $[q;0;0]$, hence, their multiplicity in $\Curv^{sm}_k$ is zero. 
\end{proof}

\begin{proof}[Proof of Theorem \ref{thm:OnHarmonicBasis}] Let us fix $\Gamma_{[\lambda]}=\Gamma^{q,p}_r$ an arbitrary $\SO(n)$-module from the previous Theorem and assume $k\leq n'/2$. Taking over the notation and slightly re-formulating the assertions from the previous proof:
$$\Hom_{\SO(n')}(\AMFP{q,q}\R^{n'}, \Res\Gamma_{[\lambda]}) = \Hom_{\SO(n')} (\bar\Upsilon^{q,0}_0 \oplus \bar\Upsilon^{q-1,0}_0, \Res \Gamma_{[\lambda]}),$$
where $1\leq \dim \Hom_{\SO(n')} (\bar\Upsilon^{q,0}_0, \Res \Gamma_{[\lambda]})\leq 2$ and $\dim \Hom_{\SO(n')} (\bar\Upsilon^{q-1,0}_0, \Res \Gamma_{[\lambda]})\leq 1$. Let us construct the basis of the space on the left-hand side.

Define $V'_{i,j}:=\AMF{i,j}\R^{n'}$, $V_\lambda := \AMF{q+r}\R^n \otimes \AMF{q}\R^n \otimes \Sym^p\R^n$. Interpreting $\SO(n')$ as the stabiliser of $\SO(n)$ which fixes $e_n\in \R^n$,  the following $\SO(n')$-equivariant map: 
\begin{eqnarray*}
\iota_{q',\lambda} : \qquad V'_{q',q'} & \to & V_\lambda \\
													v\otimes w  & \mapsto & v\wedge (e_n)^{q+r-q'}\otimes w\wedge (e_n)^{q-q'} \otimes (e_n)^{p}
\end{eqnarray*}
is injective if $q -q' + r\leq 1$ and trivial otherwise. 

Now, the map $\mu_{[q,\lambda]}:=\mu_{\lambda} \circ \pi_{\tr}\circ \iota_{q,\lambda} : V'_{q,q}\to \Res\Gamma_{[\lambda]}$ is $\SO(n')$-equivariant and its restriction to the $\SO(n')$-module $\bar\Upsilon^{q,0}_0 \subset V_{q,q}$ is not trivial. Let $v:=e_{1\ldots q}\otimes e_{1\ldots q} \in V_{q,q}$. Then $v$ fulfills all Bianchi-identities for the $\SL(n')$-module of type $[q;0;0]$, as exchanging $e_i$ from the first column with $e_j$ from the second column yields either $v$ ($i=j$) or $0$ ($i\neq j$). Hence, $\pi_{\tr}(v)\in \bar\Upsilon^{q,0}_0$ and it is straight-forward to verify that $\pi_{\tr}(v)\neq 0$. 

On the other hand, $i_{q,\lambda}(v)=:w_0$ is not a multiple of $Q:=\sum_{i=1}^n e_i^2$, since neither $e_n^2$ nor $v$ are multiples of $Q$, $v$ is not a multiple of $Q':=Q-e_n^2$, and $q\leq(n-1)/2$. Taking $\pi_{\tr}$ to be the projection on the traceless subspace with respect to $Q$, one thus obtains $\pi_{\tr}(w_0) \neq 0$. By Proposition \ref{thm:exchangeColumns}, $\mu_{\lambda}(w_0)$ is a sum of $w_0$ and several of its permutations obtained by exchanging $e_i$, $i<n$ from either the first or second column with $e_n$ from the symmetric part $e_n^p$. As the traceless part of a vector is obtained by subtracting from it certain multiples of $Q$, projecting all such permutations to trace-free spaces yields linearly independent forms. All in all, we obtain that $\mu_{[q,\lambda]}(v) \neq 0$. Hence, if $\bar\Upsilon^{q,0}_0$ is irreducible, then $\mu_{[q,\lambda]}$ spans $\Hom_{\SO(n')} (\bar\Upsilon^{q,0}_0, \Res \Gamma_{[\lambda]})$. 

As $\bar\Upsilon^{q-1,0}_0$ is always irreducible, the -- possibly, trivial -- space $\Hom_{\SO(n')} (\bar\Upsilon^{q-1,0}_0, \Res \Gamma_{[\lambda]})$ is spanned by $\mu_{[q-1,\lambda]}$. In fact, taking $v':= e_{1\ldots q-1}\otimes e_{1\ldots q-1}$ and assuming that $\iota_{q-1,\lambda}$ is not trivial, $\mu_{\lambda}(\iota_{q-1,\lambda}(v'))$ is a multiple of $\iota_{q-1,\lambda}(v')$ and from the same argument as for $\mu_{[q,\lambda]}$ follows that it contains a non-trivial traceless part. Obviously, $\mu_{[q,\lambda]}$ and $\mu_{[q-1,\lambda]}$ are linearly independent.

If $q=n'/2$, then $\bar\Upsilon^{q,0}_0 = \Upsilon^{q,0}_0 \oplus (\Upsilon^{q,0}_0)^*$ and $\dim \Hom_{\SO(n')} (\bar\Upsilon^{q,0}_0, \Res \Gamma_{[\lambda]}) = 2$. Now, the map $*_2:V'_{q,q} \to V'_{q,q}$, $(v\otimes w)\mapsto (v \otimes *w)$, where $*$ is the Hodge-operator, restricts to a non-trivial $\SO(n')$-equivariant map on $\bar\Upsilon^{q,0}_0$ which is not multiple of the identity (see \citep[p. 290]{FultonHarris}). Hence, $\mu_{[q,\lambda]}$ and $\mu^*_{[q,\lambda]}:= \mu_{[q,\lambda]}\circ *_2$ are linearly independent and span $\Hom_{\SO(n')} (\bar\Upsilon^{q,0}_0, \Res \Gamma_{[\lambda]})$.

Having the basis $\mu_{[q,\lambda]}$, $\mu_{[q,\lambda]}$ -- and $\mu^*_{[q,\lambda]}$ if $q=n'/2$ -- of $\Hom_{\SO(n')}(\AMFP{q,q}\R^{n'}, \Res \Gamma_{[\lambda]})$, let us construct an isomorphism to $(\AMFP{k,n'-k} W_y^* \otimes \Res \Gamma_{[\lambda]})^{\SO(n')}$, where $y=e_n$.

Let $V,W$ be $G$-modules for a Lie-group $G$ and $v_1,\ldots v_N$ be the basis of $V$. Any $G$-equivariant map $\mu \in \Hom_G(V,W)$ may be identified with the element $\sum_{i=1}^N v^*_i \otimes \mu(v_i) \in (V^*\otimes W)^G$. As $V'_{q,q}$ has a canonical basis $e_I \otimes e_J:= e_{i_1\ldots e_q} \otimes e_{j_1\ldots j_q}$, where $I=(1\leq i_1\leq \ldots e_q\leq n')$, we may identify $(V'_{q,q})^*$ with $V'_{q,q}$ via the map $e^*_I \mapsto e_I$ and write any $\SO(n')$-equivariant map $\mu : V'_{q,q}\to W$ as a multiple of:
$$ \bar \mu := \sum e_I \otimes e_J \otimes \mu(e_I\otimes e_J) \in (V'_{q,q}\otimes W)^{\SO(n')},$$
where the sum is over all $q$-tuples $I,J$.

Observe that the map $*_2 : V'_{i,j} \to V'_{i, n'-j}$ is an $\SO(n')$-equivariant isomorphism and so is $\nu : V'_{i,j} \to \AMF{i,j}W_y^*$ which sends $e_I\otimes e_J \mapsto \dxi{I}\otimes \dyi{J}$ for any $i$- tuple $I$ and $j$-tuple $J$. Now, $\R^{n'}\oplus\R^{n'}$ is a symplectic space with the symplectic form $Q'\in V'_{1,1}$ and the map $L^{m} : V'_{i,j}\to V'_{i+m,j+m}$ given by the $m$-fold application of the Lefschetz operator $L : V'_{i,j} \to V'_{i+1,j+1}$, $v\otimes w \mapsto (v\otimes w) \wedge Q' := \sum_{i=1}^{n'} v\wedge e_i \otimes w\wedge e_i$, is injective for $i+j \leq n'-2m$. Hence, $\nu\circ *_2 \circ L^{k-q}$ is $\SO(n')$-equivariant and injective and so is the map
\begin{eqnarray*}
\tilde \rho_{q,k,\lambda} : \Hom_{\SO(n')}(V'_{q,q},\Res\Gamma_{[\lambda]})& \to & (\AMF{k,n'-k}W_y^*\otimes \Res\Gamma_{[\lambda]})^{\SO(n')} \\
					\mu & \mapsto & \sum (\nu\circ *_2 \circ L^{k-q})(e_I\otimes e_J) \otimes \mu(e_I\otimes e_J).
\end{eqnarray*}
As $*_2 \circ L^{k-q}$ maps primitive forms to primitive forms, the restriction of $\rho_{q,k,\lambda}$ to $\AMFP{q,q}\R^{n'}$ yields the desired $\SO(n')$-equivariant isomorphism.

Note that 
$\sum L^{k-q}(e_I\otimes e_J)= \sum e_{i_1\ldots i_k}\otimes e_{j_1\ldots j_q i_{q+1}\ldots i_k},$
where the sum is over $i_1\ldots i_k, j_1\ldots j_q$. We may assume that all indexes in the sum are distinct, otherwise $L^{k-q} (e_I\otimes e_J) =0$. Hence, there is a permutation $\pi\in S_{n'}$ for each $J$ such that $(j_1\ldots j_q i_{q+1}\ldots i_k) = (\pi_1\ldots \pi_k)$. As $*e_{\pi_1\ldots \pi_k} = \sgn \pi\, e_{\pi_{k+1}\ldots \pi_n}$, one sees that $\tilde\rho_{q,k,\lambda}(\mu)$ is a multiple of
$$  \rho_{q,k,\lambda}(\mu) := \sum \sgn \pi \,\dxi{i_1\ldots i_q \pi_{q+1}\ldots \pi_k}\dyi{\pi_{k+1} \ldots \pi_n}\otimes \mu(e_{i_1\ldots i_q}\otimes e_{\pi_1\ldots \pi_q}),$$
where the sum is over $\pi\in S_{n'}$ and $i_1\ldots i_q= 1,\ldots, n'$.

All in all, the basis of $(\AMFP{k,n'-k}W_y^* \otimes \Res \Gamma_{[\lambda]})^{\SO(n')}$  consists of those elements from $\rho_{q,k,\lambda}(\mu_{[q,\lambda]})$, $\rho_{q-1,k,\lambda}(\mu_{[q-1,\lambda]})$, and $\rho_{q,k,\lambda}(\mu^{*}_{[q,\lambda]})$ which are not trivial. In particular, as $\dxi{n}|_{(0,e_n)}=\alpha$, $\dyi{n}|_{(0,e_n)}=0$ and $y_i(0,e_n)= \delta_{in}$, one has:
\begin{enumerate}
	\item If $\lambda=[q,p,0]$, $\rho_{q,k,\lambda}(\mu_{[q,\lambda]}) = \FPhi_{[k,p,q]}|_{(0,e_n)}$, $\rho_{q-1,k,\lambda}(\mu_{[q-1,\lambda]}) = \FPsi_{[k,p,q]}|_{(0,e_n)}$ and, if $q=k=n'/2$, $\rho_{q,k,\lambda}(\mu^*_{[q,\lambda]})$ is a multiple of $\FTheta_{[k,p]}|_{(0,e_n)}$;
	\item If $\lambda=[q,p,1]$, $\rho_{q,k,\lambda}(\mu_{[q,\lambda]}) = \FXi_{[k,p,q]}|_{(0,e_n)}$.
\end{enumerate}
The conditions for these forms' non-triviality may now be elaborated from the conditions for the non-triviality of $\iota_{q,\lambda}$ and Theorem \ref{thm:curvHarmonicDecomposition}. Since all $\FT_{[k,p,q]}$ are $\SO(n)$-invariant (see Remark \ref{rem:ThetaSObutNotO}, the claim now follows for all self-dual irreducible $\SO(n)$-modules $\Gamma_{[\lambda]}=\bar\Gamma_{[\lambda]}$.

If $\Gamma_{[\lambda]}$ is not self-dual, then $n=2k$ and $\lambda_k\neq 0$. Let $\lambda_k > 0$. By Remark \ref{rem:realModules}, the $\OLie(n)$-module $\bar\Gamma_{[\lambda]}$ is real. Since $*_1$ is not a multiple of the identity on  $\bar\Gamma_{[\lambda]}$, the basis of $(\Omega_p^{k,n'-k}\otimes \bar\Gamma_{[\lambda]})^{\SO(n)}$ is constituted by $\FXi_{[k,p,k]}, *_1 \FXi_{[k,p,k]}$ if $|\lambda_k|=1$ and by $\FPsi_{[k,p,k]}, *_1 \FPsi_{[k,p,k]}$ otherwise. 

In contrast, $\Gamma_{[\lambda]}$ and its dual are not always real and only complex-valued curvature measures may assume values in them. Extending $\bar\Gamma_{[\lambda]}$ to $\bar\Gamma_{[\lambda],\C}:=\bar\Gamma_{[\lambda]}\otimes\C$ by complex-linearity, one sees that $*_1$ has two eigenvalues $\pm i^m$ and the eigenspaces $E_{\pm i^m} := \{ v \mp i^m *_1 v\,|\, v\in \bar\Gamma_{[\lambda,\C]}\}$ correspond precisely to the complex $\SO(n)$-modules $\Gamma_{[\lambda]}^*$ and $\Gamma_{[\lambda]}$. This yields the claim for $m=n/2$.
\end{proof}
\begin{Remark}
\label{rem:ThetaSObutNotO}
The forms $\FT_{[k,p,q]}$, $T\in \{\Phi,\Psi,\Xi\}$, are $\SO(n)$-covariant, whereas $\Theta_{[k,p]}$ is $\OLie(n)$-covariant, as 
$$g \sum_{\pi} \sgn\pi\, y_{\pi_n} \dxi{\pi_{q+1}\ldots \pi_{k}}\dyi{\pi_{k+1}\ldots \pi_{n-1}}\otimes e_{\pi_1\ldots \pi_q} = \det g \sum_{\pi} \sgn\pi\, y_{\pi_n} \dxi{\pi_{q+1}\ldots \pi_{k}}\dyi{\pi_{k+1}\ldots \pi_{n-1}}\otimes e_{\pi_1\ldots \pi_q}$$
and $g \sum_{i=1}^n \dxi{i}\otimes e_i = \sum_{i=1}^n \dxi{i}\otimes e_i$ for all $g\in \OLie(n)$. The maps $\mu_{\pi}$ and $\pi_{\tr}$ being $\OLie(n)$-invariant do not destroy the invariances of the symmetrised differential forms. As $\nc(gK) = \det(g)\,g \nc(K)$, one has:
\begin{eqnarray*}
\CTheta_{[k,p]}(gK, gU) = \det g \int_{g (\nc(K)\cap\pi^{-1}(U))} \FTheta_{[k,p]} = \det g \int_{\nc(K)\cap \pi^{-1}(U)} g^* \FTheta_{[k,p]} = \det g\, \CTheta_{[k,p]}(K,U).
\end{eqnarray*}
In particular, $\CTheta_{[1,p]}$ is a $\Sym^p\R^3$-valued smooth translation-invariant $SO(n)$-equivariant curvature measure which is not $O(n)$-equivariant.

On the contrary, $g^* \FT_{[k,p,q]} = (\det g)\, \FT_{[k,p,q]}$ for $T\in\{\Phi,\Psi,\Xi\}$ and we obtain by the same computation as above $\CT_{[k,p,q]}(gK,gU) = \CT_{[k,p,q]}(K,U)$.
\end{Remark}

\begin{proof}[Proof of Proposition \ref{thm:formsGeometric}]
The proof requires several facts from the geometric measure theory that were also used in Section 4 of \citep{HugSchneiderTensorCurvs2013}.

Let us evaluate $\FPhi_{\otimes k,p,q}$ at the point $(x,y):=(0,e_n)$ under the assumption that the \emph{approximate tangential space} $T_{(0,e_n)} \nc(K)$ for a body $K$ has the basis $a_j := (\pd{}{x_j}, \pd{}{y_j}) \simeq ( \kappa_j b_j, \lambda_j b_j)$, $j=1,\ldots, n-1$, where $\kappa_j, \lambda_j \in [0,\infty)$ and $b_j$ is the orthonormal basis of $W:= e^{\perp}_n \subset \R^n$ with dual $b^*_j$. Then $\dxi{j} = \kappa_j\,b^*_j$ and $\dyi{j} = \lambda_j\, b^*_j$. By the skew-symmetry of the wedge-product, we see that $i_j \in \{\pi_1,\ldots, \pi_q\}$ for all $j=1,\ldots,q$, which yields at $(0,e_n)$:
$$\FPhi_{\otimes k,p,q} = (-1)^{n-1} q! \sum \sgn \pi\, \kappa_{\pi_1 \ldots \pi_k} \lambda_{\pi_{k+1} \ldots \pi_{n-1}} b^*_{\pi_1 \ldots \pi_{n-1}}\otimes (b_{\otimes \pi_1\ldots \pi_q})^{\otimes 2} \otimes y^p,$$
where the sum is over $\pi\in S_{n-1}$ and we employ the shorthand notation $\kappa_{ij} := \kappa_j \cdot \kappa_j$. Now, $\sum_\pi \sgn \pi\, e_{\otimes \pi_1\ldots \pi_q} = (q!)^{-1} \sum_\pi \sgn \pi\, e_{\pi_1\ldots \pi_q}$ and $b^*_{\pi_1\ldots \pi_{n-1}} = \sgn\pi\, \vol_W$ is just a multiple of the volume-form on $W$. Hence, 
\begin{equation}
\FPhi_{\otimes k,p,q}|_{(0,e_n)} = (-1)^{n-1} (q!)^{-1} \sum \kappa_{\pi_1 \ldots \pi_k} \lambda_{\pi_{k+1} \ldots \pi_{n-1}} \vol_W \otimes (b_{\pi_1\ldots \pi_q})^{\otimes 2} \otimes y^p,
\label{eq:PhiGeometric}
\end{equation}

We choose $\kappa_j$ and $\lambda_j$ so that $b_j$ form an orthonormal basis of $T_{(0,e_n)}$. In particular, if $b_j$ are the directions of the (generalised) principal curvatures $k_j$, then $\kappa_j = (1+k_j^2)^{-1/2}$ and $\lambda_j = k_j (1+k_j^2)^{-1/2}$ with the convention that $\kappa_j = 0$ and $\lambda_j =1$ if $k_j=\infty$. 

If $K=P$ is a polytope and $0\in F \in \mathcal{F}_s$, then there are exactly $s$ different principal curvatures $k_j$ with value 0 and exactly $(n-s-1)$ of those with value $\infty$. Hence, if $k\neq s$, then $\FPhi_{(\otimes k,p,q)}|_{(0,e_n)} = 0$. Let us now assume w.l.o.g. that $k_1=\ldots= k_k = 0$ and $k_{k+1}=\ldots \ldots = k_{n-1} = \infty$. Then $b_1,\ldots, b_s$ form the basis of $L(F)$ and $\vol_W = \vol_{L(F)} \otimes \vol_{S(F^\perp)}$, where $S(F^{\perp})$ is the unit sphere in the orthogonal complement of $L(F)$ in $\R^n$. 
$$\FPhi_{\otimes k,p,q}|_{(0,e_n)} = (-1)^{n-1} (q!)^{-1} \sum \vol_{L(F)}\otimes \vol_{S(F^{\perp})} \otimes (b_{\pi_1\ldots \pi_q})^{\otimes 2} \otimes y^p,$$
where the sum is over such $\pi \in S_{n-1}$ that $\pi_j \in \{1,\ldots, k\}$ for $j=1,\ldots, k$ and $\pi_j \in \{k+1,\ldots, n-k-1\}$ for $j=k+1,\ldots,n-k-1$. Since the term under the sum is independent of $\pi_{q+1},\ldots \pi_{n-1}$ and $\sum_{\pi \in S_k} b_{\pi_1\ldots \pi_q} = \sum_{i_1,\ldots, i_q=1}^k b_{i_1\ldots i_q}$, we see that $(b_{\pi_1\ldots \pi_q})^{\otimes 2} = Q^{\wedge q}_{L(F)}$ and obtain 
$$\FPhi_{\otimes k,p,q}|_{(0,e_n)} = (-1)^{n-1} \frac{(k-q)!(n-k-1)!}{q!} \vol_{L(F)}\otimes \vol_{S(F^{\perp})}\otimes Q^{\wedge q}_{L(F)} \otimes y^p.$$
The first identity for $\CPhi_{\otimes k,p,q}$ now follows from the definition of the Integ operator and the properties of the normal cycle for polytopes. The second one is directly implied by the first equation in the proof of \citep[Lemma 4.1]{HugSchneiderTensorCurvs2013}.
\end{proof}

\subsection{Symmetries}
\label{sec:notation}
Let us start with the following easy-to-verify identity:
\begin{equation}
\sum_{i\in\{i_1,\ldots,i_k\}}\sum_{\pi\in S_n}\sgn\pi\, e_{\pi_i}\otimes e_{\pi_{i_1}\ldots \pi_{i_k}} =  k \sum_{\pi\in S_n}\sgn\pi\, e_{\pi_{i_1}}\otimes e_{\pi_{i_1}\ldots \pi_{i_k}}.
\label{eq:mainTool}
\end{equation} 
For a $d$-partition $\mathbf{r}=(r_1,\ldots, r_d)$ of $n$, we write:
$$ \mathbf{e}_{\pi,\mathbf{r}} := e_{\pi_{s_1}\ldots\pi_{t_1}}\otimes e_{\pi_{s_2}\ldots \pi_{t_2}}\otimes \ldots \otimes e_{\pi_{s_d} \ldots \pi_{t_d}} \in \AMF{r_1}\R^n \otimes \ldots \otimes \AMF{r_d}\R^n,$$
where $t_j=\sum_{i=1}^j r_i$ and $s_j = t_{j-1} + 1$ (in particular, $s_1 = 1$ and $t_d = n$). We will refer to $e_{\pi_{s_j}\ldots \pi_{t_j}}$ as the $j$-th column or the $j$-th wedge-vector in $\mathbf{e}_{\pi,\mathbf{r}}$.

Next, define $\mathbf{e}_{\pi,\mathbf{r},i,\mathbf{k}}$, where $1\leq i\leq n$ and $\mathbf{k}\subset\{1,\ldots, d\}$, to be the vector obtained from $\mathbf{e}_{\pi,\mathbf{r}}$ by replacing the wedge-vector $e_{\pi_{s_p}\ldots \pi_{t_p}}$ with $e_{\pi_i \pi_{s_p}\ldots \pi_{t_p}}$ if $p\in\mathbf{k}$. Last, define the operation $\sigma_{pq}$ for $p\in \mathbf{k}$, $q\notin \mathbf{k}$ on $e_{\pi,\mathbf{r},i,\mathbf{k}}$ given by exchanging $e_{\pi_i}$ and $e_{\pi_{s_q}}$ in $e_{\pi_i \pi_{s_p}\ldots \pi_{t_p}}$ and $e_{\pi_{s_q}\ldots \pi_{t_q}}$.

\begin{Lemma}
\label{thm:ProtoBianchi}
Set $p\in\mathbf{k}$, write $\mathbf{k}':=\{1,\ldots, d\}\setminus \mathbf{k}$, and assume that $\mathbf{k},\mathbf{k}'$ are non-empty. Then:
\begin{equation}\label{eq:ProtoBianchi}
(r_p+1)\sum_{i,\pi}\sgn\pi\,\mathbf{e}_{\pi,\mathbf{r},i,\mathbf{k}} = \sum_{q\in \mathbf{k}'} r_q  \sum_{i,\pi}\sgn\pi\, \sigma_{pq}(\mathbf{e}_{\pi,\mathbf{r},i,\mathbf{k}}),
\end{equation}
where the sum is over $\pi\in S_n$ and $i=1,\ldots,n$.
\end{Lemma}
\begin{proof} We may re-order the wedge-vectors in $\mathbf{e}^{\pi}_{\mathbf{r},i,\mathbf{k}}$ and assume $\mathbf{k}=(1,\ldots, d-u)$, $\mathbf{k}'=(d-u+1,\ldots, d)$ for $1 < u < d$, and $p=1$. The proof will now be carried out inductively over $|\mathbf{k}'|=u$. For the sake of brevity, we omit the subscript $\mathbf{k}$ in $\mathbf{e}_{\pi,\mathbf{r},i,\mathbf{k}}$ in the proof.

Let $|\mathbf{k}'|=1$ and, hence, $\mathbf{k}'=(d)$. As $e_{ii} = 0$, we have:
$$\sum_{i,\pi}\sgn\pi\,\mathbf{e}_{\pi,\mathbf{r},i} = \sum_{i=s_d}^{t_d} \sum_{\pi}\sgn\pi\,\mathbf{e}_{\pi,\mathbf{r},i} \overset{\pref{eq:mainTool}}{=}  r_q \sum_{\pi} \sgn\pi\,\mathbf{e}_{\pi,\mathbf{r},s_d} =: I.$$
All wedge-vectors of $\mathbf{e}_{\pi,\mathbf{r},s_d}$ begin with the vector $e_{\pi_{s_d}}$, hence $\sigma_{1 d} (\mathbf{e}_{\pi,\mathbf{r},s_d}) = \mathbf{e}_{\pi,\mathbf{r},s_d}$ and:
$$ I =r_q \sum_{\pi} \sgn\pi\, \sigma_{1d}(\mathbf{e}_{\pi,\mathbf{r},s_d}) \overset{\pref{eq:mainTool}}{=} \frac{r_q}{r_1+1} \sum_{i\in{s_1,\ldots,t_1,s_d}} \sum_{\pi} \sgn\pi\, \sigma_{1d}(\mathbf{e}_{\pi,\mathbf{r},i})$$
We now add $0=\sgn\pi\, \sigma_{1d}(\mathbf{e}_{\pi,\mathbf{r},i})$ for $s_2\leq i\leq t_d$, $i\neq s_d$ and conclude the proof for $|\mathbf{k}'| =1$.

Assuming the claim's validity for all $|\mathbf{k}'|=u-1$, the proof for $|\mathbf{k'}| = u$ works as follows. We start by splitting the sum:

$$ \sum_{i,\pi} \sgn\pi\,\mathbf{e}_{\pi,\mathbf{r},i} = \sum_{\pi} \sgn\pi\left(\sum_{i=s_1}^{t_{d-1}}  \mathbf{e}_{\pi,\mathbf{r},i} + \sum_{i=s_d}^{t_d}  \mathbf{e}_{\pi,\mathbf{r},i}\right) =: A+B.$$
Now, $\mathbf{e}_{\pi,\mathbf{r},i} = \mathbf{e}_{\pi,\mathbf{r}',i}\otimes e_{\pi_{s_d}\ldots \pi_{t_d}}$, where $\mathbf{r}'=\mathbf{r}\setminus \{r_d\}=(r_1,\ldots, r_{d-1})$. As $|\{1,\ldots d-1\}\setminus \mathbf{k}|=t-1$. We may apply the Lemma on $\mathbf{e}_{\pi,\mathbf{r}',i}$ in $A$, observe that $\sigma_{1q}(\mathbf{e}_{\pi,\mathbf{r}',i})\otimes e_{\pi_{s_d}\ldots \pi_{t_d}} = \sigma_{1q}(\mathbf{e}_{\pi,\mathbf{r},i})$ for $q\leq d-1$, and add $0=\sum_{i=s_d}^{t_d} \sigma_{1q}(\mathbf{e}_{\pi,\mathbf{r},i}) - \sum_{i=s_d}^{t_d} \sigma_{1q}(\mathbf{e}_{\pi,\mathbf{r},i})$ to obtain:

$$ A = \sum_{q\in \mathbf{k}'\setminus\{d\}}  \frac{r_q}{r_1+1} \left(\sum_{i,\pi} \sgn\pi\,\sigma_{1q}(\mathbf{e}_{\pi,\mathbf{r},i}) - \sum_{\pi} \sgn\pi \sum_{i=s_d}^{n} \sigma_{1q}(\mathbf{e}_{\pi,\mathbf{r},i})\right).$$
The second summand may be re-written for any $q \in \mathbf{k}'\setminus\{d\}$:

$$ \sum_{\pi}\sgn\pi \sum_{i=s_d}^{t_d} \sigma_{1q}(\mathbf{e}_{\pi,\mathbf{r},i}) \overset{\pref{eq:mainTool}}{=} r_d \sum_{\pi} \sgn\pi\, \sigma_{1q}(\mathbf{e}_{\pi,\mathbf{r},s_d}) \overset{\pref{eq:mainTool}}{=} -\frac{r_d}{r_q}  \sum_{\pi}\sgn\pi\sum_{i=s_q}^{t_q} \sigma_{1q}(\mathbf{e}_{\pi,\mathbf{r},i}),  $$
since $\sum_{\pi} \sgn\pi e_{\pi_{s_q}} \otimes (e_{\pi_{s_d}})^{\otimes d-u} = -\sum_{\pi} \sgn\pi e_{\pi_{s_d}} \otimes (e_{\pi_{s_q}})^{\otimes d-u}$. As in the case $|\mathbf{k}'|=1$,
$$ B = \frac{r_d}{r_1+1} \sum_{\pi} \sgn\pi\left(\sum_{i=1}^{t_{d-u}} \sigma_{1d}(\mathbf{e}_{\pi,\mathbf{r},i}) + \sum_{i=s_d+1}^{t_d} \sigma_{1d}(\mathbf{e}_{\pi,\mathbf{r},i})\right),$$
which concludes the proof for all $p,n,\mathbf{r},\mathbf{k}$.
\end{proof}

To prove Theorem \ref{thm:BernigLemma}, we need a finer control over the symmetrisation of forms. We write $\FT^{\pi}_{\otimes k,p,q}:= \FT_{\otimes k,p,q} \cdot \pi$ for the forms obtained by permuting its tensor part by some permutation $\pi \in S_{|\lambda|}$ of the Young-diagram $\lambda=[q;p;r]$ as in \pref{eq:YoungForms}. More generally, we write $\FT^{d}_{\otimes k,p,q}:= \FT_{\otimes k,p,q} \cdot d$ for any symmetrisation by an element $d$ of the group algebra $\C S_{|\lambda|}$. For the sake of brevity, we will write $\pi$ instead of $e_\pi$ for the basis elements of $\C S_{|\lambda|}$. 

There are several distinguished permutations. We write $(i_a\,j_b) \in S_{|\lambda|}$ for the transposition which exchanges the $a$-th box in the $i$-th column with the $b$-th box in the $j$-th column and $\sigma_{\ell}:=\prod_{j=1}^\ell (1_j\,2_j)$ for the permutation which exchanges the first $\ell\leq q$ boxes in the first column with the same number of boxes in the second column. More generally, define $\sigma_{\mathbf{r}} := \prod_{j\in\mathbf{r}} (1_j\,2_j)$ for any subset $\mathbf{r}\in \{1,\ldots, q\}$ and $d_{\lambda}:=\id+\sigma_q\in \C S_{|\lambda|}$. 

As all eligible Young-diagrams $[q;p;r]$ have at most one box in any column starting with the third, we write $j$ instead of $j_1$ for any $j\geq 3$. Let $R'_\lambda$ be the group of permutations generated by transpositions $(i\, j)$, $i,j\geq 3$, and $h_{\lambda}:= \sum_{\pi\in R'_\lambda} \pi$ and define the following symmetrised forms:
\begin{equation} \label{eq:SLnAndSOnCurvs}
\FT^{\pi}_{k,p,q} := \FT^{\pi\cdot b_\lambda}_{\otimes k,p,q}, \qquad \FT^{\pi}_{(k,p,q)} := \FT^{\pi\cdot h_{\lambda} \cdot b_\lambda}_{\otimes k,p,q}, \qquad \FT^{\pi}_{\{k,p,q\}} := \FT^{\pi\cdot a_\lambda\cdot b_\lambda}_{\otimes k,p,q},
\end{equation}
where $a_\lambda, b_\lambda$ are as in eq. \pref{eq:YoungSymmetrisers}. They assume values in $\AMF{\lambda'}\R^n=\AMF{q+r, q}\R^n \otimes (\R^n)^{\otimes p}$, $\AMF{q+r, q}\R^n \otimes\Sym^p \R^n$, and $\Gamma_{\lambda}$, respectively.  Note that $\FT_{k,p,q}$ satisfy the following lower-rank relations:
\begin{equation}
\FPhi_{k,1,0} = \FXi_{k,0,0} \qquad\hbox{and}\qquad \FPsi_{k,p,1} = \FXi_{k,p+1,0} = \FPhi_{k,p+2,0}.
\label{eq:FamilyRelations}
\end{equation}
We use the same notation for the symmetrisations of the curvature measures $T^{\pi}_{\otimes k,p,q}$. 
\begin{Ex} \label{ex:formSymmetrisation}As $y^{p}\cdot h_\lambda = p!\,y^p$, one has $ \FT_{(k,p,q)} = p!\, \FT_{k,p,q}$ for $\FT\in \{\FPhi, \FXi,\FPsi\}$. Similarly:
\begin{equation*}
\begin{aligned}
\FPhi^{(1_1\,3)}_{(k,p,q)} = & C_{p-1} \sum_{\pi,i}\sgn\pi\, y_{\pi_n}\dxi{i_1\ldots i_q \pi_{q+1}\ldots \pi_k}\wedge \dyi{\pi_{k+1}\ldots \pi_{n-1}}\otimes e_{y i_2\ldots i_q}\otimes e_{\pi_1\ldots \pi_q}\otimes e_{i_1}  y^{p-1}\\
\FPhi^{(1_1\,3)\cdot \sigma_q}_{(k,p,q)} = & C_{p-1} \sum_{\pi,i}\sgn\pi\, y_{\pi_n}\dxi{i_1\ldots i_q \pi_{q+1}\ldots \pi_k}\wedge \dyi{\pi_{k+1}\ldots \pi_{n-1}}\otimes e_{\pi_1\ldots \pi_q}\otimes e_{y i_2\ldots i_q}\otimes e_{i_1}  y^{p-1}\\
\FPhi^{(1_1\,3)(2_1\,4)}_{(k,p,q)} = & C_{p-2} \sum_{\pi,i}\sgn\pi\, y_{\pi_n}\dxi{i_1\ldots i_q \pi_{q+1}\ldots \pi_k}\wedge \dyi{\pi_{k+1}\ldots \pi_{n-1}}\otimes e_{y i_2\ldots i_q}\otimes e_{y\pi_2\ldots \pi_q}\otimes e_{i_1} e_{\pi_1}  y^{p-2}, 
\end{aligned}
\end{equation*}
where the sums are as in \pref{eq:FormDefs} and $C_p = (-1)^{n-1} p!$. 
\end{Ex}

\begin{Prop} \label{thm:exchangeColumns}
For any $\FT\in\{\FPhi,\FXi,\FPsi\}$ and $\mathbf{r}\subset\{1,\ldots, q\}$, one has:
\begin{equation}
\FT^{\sigma_\mathbf{r}}_{k,p,q} \equiv  \binom{q'}{\ell'}^{-1} \FT_{k,p,q} \quad \mod \{d\alpha\} \label{eq:exchangeColumns},
\end{equation}
where $q'= q-1$ if $\FT=\FPsi$ and $q$ otherwise, and $\ell'\leq q'$ is the number of transpositions in $\sigma_\mathbf{r}$ which exchange $e_{i_a}$ with $e_{\pi_a}$.
Furthermore, one has:
\begin{eqnarray} \label{eq:SLnBase}
\begin{aligned}
\CPsi_{\{k,p,q\}}  =\, & 2 q\, \CPsi_{(k,p,q)}, \quad \CXi_{\{k,p,q\}} = (q+1) \CXi_{(k,p,q)} + qp\,\left(\CXi^{(2_1\,3)}_{(k,p,q)} - \frac{q-1}{2} \CXi^{(1_1\,2_1\,3)}_{(k,p,q)}\right), \\
\quad \CPhi_{\{k,p,q\}}&=  (q+1) \Phi_{(k,p,q)} + q p \left(\CPhi^{(2_1\,1_1\,3)\cdot d_\lambda}_{(k,p,q)} + \CPhi^{(2_1\,3)\cdot d_\lambda}_{(k,p,q)} + (p-1)\Phi^{(1_1\,3)(2_1\,4)}_{(k,p,q)}\right).
\end{aligned}
\end{eqnarray}
\end{Prop}

\begin{proof}
As $\FT^{(1_a\,2_a)}_{k,p,q} = \FT^{(1_{b}\, 2_{b})}_{k,p,q}$ for all $a,b\leq q$, we may assume $\mathbf{r}=(1,\ldots, \ell)$ and $\sigma_{\mathbf{r}}=\sigma_\ell$. By the $\SO(n)$-covariance of the forms, it suffices to show the claim for at the point $(0,e_n)$. As the above permutations exchange the boxes contained in the first two columns and $\FPsi_{\otimes k,p,q}^{(1_q\, 2_q)}= \FPsi_{\otimes k,p,q}$, it suffices to prove \pref{eq:exchangeColumns} for $Z:=\FPhi_{\otimes k,0,q}|_{(0,e_1)}$ with $\ell'=\ell$ and $q'=q$. We do this by induction over $\ell$. The case $\ell=0$ is trivial. Now assume that the claim is valid for $\ell-1$. Set
\begin{eqnarray*}
Y_1 & = & \sum \sgn\pi\, \dxi{i_\ell \pi_{q+1}\ldots \pi_{k}}\otimes \dyi{\pi_{k+1}\ldots \pi_{n-1}}\otimes e_{\pi_1\ldots \pi_{\ell-1} i_\ell} \otimes e_{\pi_{\ell}\ldots \pi_{q}}
\end{eqnarray*}
and let $Y_2$ be the element obtained by exchanging $e_{i_\ell}$ and $e_{\pi_{\ell}}$. Then, by Lemma \ref{thm:ProtoBianchi}, $\ell Y_1 \equiv (q-l+1) Y_2$ mod $d\alpha$. Furthermore, $Z^{\sigma_{\ell-1}\cdot b_\lambda}$ and $Z^{\sigma_\ell \cdot b_\lambda}$ are the images of $Y_1$ and $Y_2$ under the \emph{injective} map which wedges $q-\ell$ copies of $Q':=\sum \dxi{i}\otimes e_i$ with the first and the third columns and $\ell-1$ copies of $Q'$ with the first and the fourth column. We conclude:
\begin{eqnarray*}
Z^{\sigma_\ell \cdot b_\lambda} & \equiv & \frac{l}{q-l+1} Z^{\sigma_{\ell-1} \cdot b_{\lambda}} \equiv \binom{q}{\ell}^{-1} Z^{b_\lambda}\quad \mod d\alpha.
\end{eqnarray*}

Let us analyse the structure of $c_\lambda$ for $\lambda=[q;p;r]$. It is clear that $a_\lambda = \prod_{j=1}^q a_j$, where $a_j$ is the sum over the elements from $S_{|\lambda|}$ which preserve the $j$-th row. 

Setting $d_j = \id + (1_j\,2_j)$, we see that $a_j = d_j$ if $j\geq 2$. On the contrary, the subgroup of $S_{|\lambda|}$ which preserves the first row is isomorphic to $S_{p+2}$, as there are $p+2$ boxes in the first row. Writing $S_{p+2} \simeq R'' \cdot R'_\lambda$, where $R''$ is the set of representatives of all $(p+1)(p+2)$ right cosets in $S_{p+2}/S_p$ and setting $ R'' := \{\id, (2_1\, b)\}\times  \{\id, (1_1\,2_1), (1_1\, b)\},$ where $3\leq b\leq p+2$ in both subsets:
\begin{equation}
a_1= \left(\id+\sum (2_1\, b)\right) \cdot \left(d_1 + \sum (1_1\, b) \right) \cdot h_\lambda,
\label{eq:a1rewrite}
\end{equation}
where the sums are over $b=3,\ldots, p+2$ and $a'_1:= \sum_{h\in R'} h$. As $(2_1\, b)(1_1\, b)=(1_1\, b)(1_1\,2_1)$, the first two terms can be re-written as $[\id + \sum ((1_1\, b) + (2_1\, b))]\cdot d_1 + \sum_{b\neq b'} (1_1\, b)(2_1\, b')$, where $b,b'$ run from $3$ to $p+2$. As $h_\lambda$ symmetrises all columns beginning with the third, we have for all $i,j\in\{1,2\}$, $i\neq j$ and $3\leq b \leq p+2$:
\begin{eqnarray*}
\begin{matrix*}[l]
\FT^{(j_1\, b)\cdot h_\lambda}_{\otimes k,q,p} =  \FT^{(j_1\, 3)\cdot h_\lambda}_{\otimes k,q,p}, & \quad\FT^{(i_1\, b)(j_1 b')\cdot h_\lambda}_{\otimes k,q,p} = \FT^{(i_1\, b')(j_1\, b)\cdot h_\lambda}_{\otimes k,q,p}, & \quad\FT^{(i_1\,j_1\, b)\cdot h_\lambda}_{\otimes k,q,p} =  \FT^{(i_1\,j_1\, 3)\cdot h_\lambda}_{\otimes k,q,p} 
\end{matrix*}
\end{eqnarray*}
and $a_1= \left(\id+ p (1_1\, 3) + p (2_1\, 3) + \frac{p(p-1)}{2} (1_1\, 3)(2_1\, 4) \right) \cdot d_1 \cdot h_\lambda$. As $h_\lambda$ and $d_j$ commute, we obtain:
\begin{equation}
c_\lambda = a_\lambda\cdot b_\lambda = \left(\id+ p (1_1\, 3) + p (2_1\, 3) + \frac{p(p-1)}{2} (1_1\, 3)(2_1\, 4) \right)\cdot d'_\lambda \cdot h_\lambda \cdot b_\lambda,
\label{eq:alambda}
\end{equation} 
where $d'_\lambda := \prod_{j=1}^q d_j = \sum_{|\mathbf{r}| \leq q} \sigma_{\mathbf{r}}$. Applied on $\CPhi_{\otimes k,p,q}$, this yields:
$$\CPhi_{\{k,p,q\}}= \CPhi^{c_\lambda}_{\otimes k,p,q} = \CPhi^{d'_\lambda}_{(k,p,q)}+ p\, \CPhi_{(k,p,q)}^{(1_1\, 3)\cdot d'_\lambda} + p\, \Phi_{(k,p,q)}^{(2_1\, 3)\cdot d'_\lambda} + \frac{p(p-1)}{2} \CPhi_{(k,p,q)}^{(1_1\, 3)(2_1\, 4)\cdot d'_\lambda}.$$
All we need to do is to compute $\Phi^{\pi\cdot d'_\lambda}_{\otimes k,p,q}$ for four different permutations $\pi$. By eq. \pref{eq:exchangeColumns}:
\begin{eqnarray*}
\CPhi^{d'_\lambda}_{(k,p,q)}\mkern-10mu & = & \sum_{\ell =0}^{q}\sum_{|\mathbf{r}|=\ell} \CPhi^{\sigma_\mathbf{r}}_{(k,p,q)}= \sum_{\ell =0}^{q} \sum_{|\mathbf{r}|=\ell} \binom{q}{\ell}^{-1} \CPhi_{(k,p,q)} = \sum_{\ell =0}^{q} \binom{q}{\ell} \binom{q}{\ell}^{-1} \CPhi_{(k,p,q)}  = (q+1)\,\CPhi_{(k,p,q)}.
\end{eqnarray*}
Similarly, one obtains $\Phi_{(k,p,q)}^{(1_1\, 3)(2_1\, 4)\cdot d'_\lambda} = 2q\, \Phi_{(k,p,q)}^{(1_1\, 3)(2_1\, 4)}$. To compute the remaining two summands, we re-write $d'_\lambda$ as follows. Set $d(a):= \prod_{j=1,j\neq a}^q d_j$ for $a\leq q$ and $\mathbf{r}^\perp := \{1,\ldots,q\}\setminus \mathbf{r}$. Observing that $\sigma_{\mathbf{r}^\perp} = \sigma_q \circ \sigma_\mathbf{r} = \sigma_{\mathbf{r}}\circ \sigma_q$, we have:
$$d'_\lambda = (1+\sigma_q) + \sum_{\ell =1}^{q-1} \sum_{|\mathbf{r}|=\ell} \sigma_{\mathbf{r}} = d_\lambda + \sum_{\ell=1}^{q-1} \sum_{|\mathbf{r}|=\ell, a\notin \mathbf{r}} \sigma_{\mathbf{r}} + \sigma_{\mathbf{r}^\perp}  = \sum_{l=0}^{q-1}\sum_{|\mathbf{r}|=\ell, a\notin \mathbf{r}} \sigma_{\mathbf{r}} \cdot d_\lambda = d(a) \cdot d_\lambda.$$
Then one sees $\CPhi_{(k,p,q)}^{(1_1\, 3)\cdot d(1)} = \frac{q+1}{2}\CPhi_{(k,p,q)}^{(1_1\, 3)}$ and $\CPhi_{(k,p,q)}^{(2_1\, 3)\cdot d(1)} = q \CPhi_{(k,p,q)}^{(2_1\, 3)} - \frac{q-1}{2} \CPhi_{(k,p,q)}^{(1_1\, 3)\cdot \sigma_q}$. As $d_\lambda = \id + \sigma_q$ on $\CPhi_{(k,p,q)}^{(1_1\, 3)\cdot d(1)} + \CPhi_{(k,p,q)}^{(2_1\, 3)\cdot d(1)}$ and $\CPhi_{k,p,q}^{(1_1\, 3)}= q  \CPhi_{k,p,q}^{(2_1\,1_1\, 3)}$, we obtain the claim for $\FPhi_{\{k,p,q\}}$. 

The computation is simpler for $\Xi,\Psi$. As there may be at most one $y$ in each column, one has: 
$$ \FXi^{(1_1 3)(2_1 4)}_{(k,p,q)} = \FXi^{(1_1 3)}_{(k,p,q)} = \FPsi^{(1_1 3)}_{(k,p,q)} =\FPsi^{(2_1 3)}_{(k,p,q)} = \FPsi^{(1_1 3)(2_1 4)}_{(k,p,q)} = 0.$$
The remaining terms are computed as above.
\end{proof}

\subsection{Globalisation}
\label{sec:Globalisation} 

\begin{proof}[Proof of Theorem \ref{thm:BernigLemma}]
To prove \pref{eq:GlobPhiPsi}, consider the $\SO(n)$-equivariant section:
$$h_{k,n} = \frac{1}{n-k-1}\sum_{j=1}^n \pd{}{y^j} \otimes y \otimes e_j \in \Gamma(T\R^n\otimes (\R^n)^{\otimes 2})$$
and set $\tilde E_{k,p,q+1}:=-\iota_{h_{k,n}} \FPhi_{(k,p,q)}$, where $\pd{}{y^j}$ is contracted with the differential form and $y \otimes e_j$ is wedged with the first two columns of its tensor-part. Then:
$$\tilde E_{k,p,q+1} = (-1)^{k+1}C_p\sum \sgn{\pi}\ y_{\pi_{n}}\, \dxi{i_1\ldots i_q \pi_{q+1}\ldots \pi_k}\wedge\dyi{\pi_{k+2}\ldots\pi_{n-1}} \otimes\,e_{i_1\ldots i_{q}y}\otimes e_{\pi_{1}\ldots\pi_{q}\pi_{k+1}} \otimes y^{p},$$
where the sum is over $i_1,\ldots, i_q = 1,\ldots, n$ and $\pi\in S_n$ and $C_p$ is as in Example \ref{ex:formSymmetrisation}. Then:
\begin{eqnarray*}
d\tilde E_{k,p,q+1}\mkern-12mu& = &\mkern-10mu -C_p\Big[\sum \sgn \pi\, \dxi{i_1\ldots i_q \pi_{q+1}\ldots \pi_{k}}\wedge\dyi{\pi_{n}\pi_{k+2}\ldots \pi_{n-1}} \otimes \,e_{i_1\ldots i_{q}y}\otimes e_{\pi_{1}\ldots\pi_{q}\pi_{k+1}} \otimes y^{p} \\
																&  & + \sum_{j=1}^n \sgn \pi\, y_{\pi_{n}}\dxi{i_1\ldots i_q \pi_{q+1}\ldots \pi_{k}}\wedge\dyi{j\pi_{k+2}\ldots \pi_{n-1}} \otimes \,e_{i_1\ldots i_{q}y}\otimes e_{\pi_{1}\ldots\pi_{q}\pi_{k+1}} \otimes y^{p-1} e_j  \\
																&  & + \sum_{j=1}^n \sgn \pi\, y_{\pi_{n}}\dxi{i_1\ldots i_q \pi_{q+1}\ldots \pi_{k}}\wedge\dyi{j\pi_{k+2}\ldots \pi_{n-1}} \otimes \,e_{i_1\ldots i_{q}j}\otimes e_{\pi_{1}\ldots\pi_{q}\pi_{k+1}} \otimes y^{p}\Big].
\end{eqnarray*}
 After having computed the exterior derivative, we may restrict the forms to $(x,y)=(0,e_n)$. Lemma \ref{thm:ProtoBianchi} and Example \ref{ex:formSymmetrisation} yield:
\begin{eqnarray*}
 d\tilde E_{k,p,q+1} \equiv  (q+1) \FPsi_{(k,p,q+1)} + \frac{p(k-q)}{n-k-1} \FPhi^{(1_1\,3)}_{(k,p,q+1)} + \frac{k-q}{n-k-1} \FPhi_{(k,p,q+1)}
\end{eqnarray*}
or, after multiplying by $(n-k-1)$ and replacing $q$ with $q-1$,
\begin{equation*}
(n-k-1)\, d\tilde E_{k,p,q} \equiv q(n-k+1) \FPsi_{(k,p,q)} + (k-q+1)\, \FPhi_{(k,p,q)} + p(k-q+1)\, \FPhi^{(1_1\,3)}_{(k,p-1,q)},
\end{equation*}
where the equality again holds modulo multiples of $\alpha$, $d\alpha$ and the forms whose tensor-parts are multiples of $Q$.

Let us now apply $\integ \otimes \pi_{\tr}\circ \mu_\lambda$ on both sides of the above equation, where $\lambda=[q;p;0]$. Recall that $\integ$ eliminates all exact forms and multiples of $\alpha,d\alpha$. Thus, one has:
\begin{equation*}
0 \equiv q(n-k+1) \VPsi_{(k,p,q)} + (k-q+1)\, \VPhi_{(k,p,q)} + p(k-q+1)\, \VPhi^{(1_1\,3)}_{(k,p-1,q)},
\end{equation*}
where the equality now holds only up to the forms whose tensor-parts are multiples of $Q$.

Applying $\mu_\lambda$ on the tensor-part, we have similarly to Example \ref{ex:formSymmetrisation}:
$$ \FPhi_{(k,p,q)}\cdot c_\lambda = \FPhi_{\otimes k,p,q} \cdot h_\lambda\cdot b_\lambda \cdot a_\lambda\cdot b_\lambda = p! (q!)^2\, \FPhi_{\otimes k,p,q}\cdot a_\lambda\cdot b_\lambda  = p!\, q!^2\, \FPhi_{\{k,p,q\}}.$$
One shows similarly to the proof of equation \pref{eq:a1rewrite} that $b_\lambda := b_{1,q}\cdot b_{2,q}$ with $b_{i,j}$, $i=1,2$, defined recursively by $\textstyle b_{i,j}=b_{i,j-1}\cdot b'_{i,j,a}$, where $b'_{i,j,a}:= \id - \sum_{r=1,r\neq a}^{j} (r_i\,q_i)$ for any $a\in \{1,\ldots j\}$, and $b_{i,1} = \id$. Using this identity, one obtains:
\begin{eqnarray*}
\FPsi_{(k,p,q)}\cdot c_\lambda & = & p!  (q-1)!^2\, \FPsi_{\otimes k,p,q}\cdot b'_{1,q,q}\cdot b'_{2,q,q} \cdot c_\lambda = p!\, q!^2\,\FPsi_{\{k,p,q\}}, \\
\FPhi^{(1_1\,3)}_{(k,p,q)}\cdot c_\lambda & = & p! q! (q-1)! \, \FPhi^{(1_1\,3)}_{\otimes k,p,q}\cdot b'_{1,q,1} \cdot c_\lambda = p!\, q!^2\,\FPhi^{(1_1\,3)}_{\{k,p,q\}}
\end{eqnarray*}
By Proposition \ref{thm:BianchiRelationsSAV}, one sees $q\,\FPhi^{(1_1\,3)}_{\{k,p,q\}} = \FPhi_{\{k,p,q\}}$. Applying $\pi_{\tr}$ which eliminates all forms whose tensors are multiples of $Q$, we obtain by eq. \pref{eq:HarmonicForm} the identity in \pref{eq:GlobPhiPsi}.

The cases \pref{eq:GlobXi} and \pref{eq:GlobPsiMargin} follow immediately from the Alesker-Bernig-Schuster decomposition of $\Val_k$ and \ref{thm:curvHarmonicDecomposition}, as the corresponding curvature measures $\CXi_{[k,p,q]}$ and $\CPsi_{[k,p,k+1]}$ assume values in $\SO(n)$-modules which occur in $\Curv^{sm}_k$ but are missing in $\Val_k$. To prove $\glob \CTheta_{[p]} = 0$ observe that $\glob \CPhi_{[k,p,k]}$ is a non-trivial $\OLie(n)$-covariant valuation with values in the same module $\Gamma^{k,p}_0$ as $\glob \CTheta_{[p]}$. As $\dim \Val^{\SO(n)}_{k,\Gamma^{k,p}_0} = 1$, all $\Gamma^{k,p}_0$-valued valuations of degree $k$ are $\OLie(n)$-invariant in contrast to $\glob \CTheta_{[p]}$ which is $\SO(n)$- but not $\OLie(n)$-covariant by Remark \ref{rem:ThetaSObutNotO}.
\end{proof}

\begin{Prop}\label{thm:smoothness}
Continuous $\Gamma$-valued $\SO(n)$-covariant translation-invariant valuations are smooth for any finite-dimensional $\SO(n)$-module $\Gamma$.
\end{Prop}
\begin{proof}
Let $\phi$ be a $\Gamma$-valued valuation satisfying the conditions in the claim. Since $\Val^{sm}$ lies dense in $\Val$, we may find a sequence $\phi_i$ of smooth $\Gamma$-valued translation-invariant valuations which converges componentwise to $\phi$. Define the map $A$ for any translation-invariant $\Gamma$-valued valuation $\tau$: 
\begin{eqnarray*}
(A\tau)(K) := \int_{\SO(n)} g^{-1} \tau(gK)\, dg.
\end{eqnarray*}
If $\tau$ is smooth, then so is $A\tau(K)$. Furthermore, for any $h\in \SO(n)$, one has
$$
A\tau(hK) =\int_{\SO(n)} g^{-1} \tau(g h K)\, dg \overset{\tilde g:=gh}{=} \int_{\SO(n)} (\tilde g h^{-1})^{-1} \phi(\tilde g K)\, d\tilde g = h (A\tau(K)),
$$
i.e. $A\tau$ is also $\SO(n)$-covariant. Applying $A$ to both the sequence $\phi_i$ and $\phi$, one obtains a sequence $A\phi_i$ of smooth $\SO(n)$-covariant translation-invariant valuations converging to $A\phi = \phi$. We have seen in the previous Sections that the space of smooth $\Gamma$-valued $\SO(n)$-covariant translation-invariant valuations is finite-dimensional and, thus, closed. Hence, $\phi =\lim_{i} A\phi_i$ is also smooth and the result follows.
\end{proof}

\begin{proof}[Proof of Proposition \ref{thm:TValBasis}]
We know that $\glob : \TCurv^{sm, \SO(n)}_{k,\Gamma_{[\lambda]}} \to \TVal_{\Gamma_{[\lambda]}}^{k,\SO(n)}$ is surjective. Let us work out its kernel. The elements $\xi^n_{[k,p,q]}$ and $\theta^n_{[p]}$ belong to the kernel by \pref{eq:GlobXi} and the elements $\psi^n_{[k,p,q]}$ either lie in the kernel by \pref{eq:GlobPsiMargin} or $\glob \psi^n_{[k,p,q]} = C_{n,k,p,q} \glob \phi^n_{[k,p',q']}$ for some constant $C_{n,k,p,q}$ and some $p'$ and $q'$ by \pref{eq:FamilyRelations} or \pref{eq:GlobPhiPsi}. The only exception is $\psi^n_{[k,p,k+1]}$ for $\frac{n-1}{2}< k \leq n-1$, as none of the relations apply to them.

The coefficients of all linearly independent $\Gamma_{[\lambda]}$-valued valuations $\VT_{[k,p,q]}$ span the isotypical component $\Gamma_{[\lambda]}$ in the space $\Val^{f}_{k}$ of the so-called $\SO(n)$-\emph{finite vectors} in $\Val_k$.  We refer to \citep[Section 3.2]{Sepanski} for the details on $G$-finite vectors in infinite-dimensional representations. As, by Alesker's Irreducibility Theorem, $\Val^{f}$ lies dense in $\Val^{sm}$, we obtain the claim.
\end{proof}

\addcontentsline{toc}{section}{References}

\bibliography{references}

\end{document}